\documentclass[11pt]{amsart}
\usepackage{mathtools,amssymb} 
\usepackage[T1]{fontenc} 
\usepackage{lmodern} 
\usepackage{graphicx} 
\usepackage[font=small,labelfont=bf]{caption} 
\usepackage{xcolor} 
\usepackage{hyperref} 
\hypersetup{
    colorlinks,
    allcolors={blue!50!black},
    linkcolor={red!50!black},
}

\renewcommand*{\subsubsection}[1]{\medskip\noindent\textbf{#1.}}

\newcommand{\R}{\mathbb{R}} 
\newcommand{\T}{\mathbb{T}} 
\newcommand{\N}{\mathbb{N}} 
\newcommand{\Z}{\mathbb{Z}} 
\newcommand{\RP}{\mathbb{RP}} 
\DeclareMathOperator{\vol}{vol} 
\DeclareMathOperator{\ch}{ch} 
\DeclareMathOperator{\dias}{dias} 
\DeclareMathOperator{\area}{area}
\DeclareMathOperator{\vertex}{vertex}
\DeclareMathOperator{\weight}{weight}
\newcommand{\D}{\mathcal{D}} 
\newcommand{\Sph}{\mathcal{S}} 

\newtheorem{theorem}{Theorem}[section]
\newtheorem{corollary}[theorem]{Corollary}
\newtheorem{lemma}[theorem]{Lemma}
\newtheorem{proposition}[theorem]{Proposition}
\newtheorem*{mahlerconjecture}{Mahler conjecture}
\newtheorem{innerconjecture}{Conjecture}
\newenvironment{conjecture}[1]
    {\renewcommand\theinnerconjecture{#1}\innerconjecture}
    {\endinnerconjecture}
\newtheorem{innertheorem}{Theorem}
\newenvironment{maintheorem}[1]
    {\renewcommand\theinnertheorem{#1}\innertheorem}
    {\endinnertheorem}
\theoremstyle{definition}
\newtheorem{definition}[theorem]{Definition}
\newtheorem*{definition*}{Definition}
\newenvironment{claim}
    {\smallskip\noindent\textbf{Claim.}\itshape}
    {}
\newenvironment{claimproof}
    {\smallskip$\rhd$}
    {$\lhd$\smallskip}
\newcommand{\term}[1]{\emph{#1}}
\newcommand{\assertion}[1]{\textit{#1}}

\title{Isosystolic inequalities for optical hypersurfaces}

\author{J.C. \'Alvarez Paiva}

\address{J.C. \'Alvarez Paiva, Laboratoire Paul Painlev\'e, Bat. M2, Universit\'e des Sciences et Technologies, 59 655 Villeneuve d'Ascq, France.}

\email{juan-carlos.alvarez-paiva@math.univ-lille1.fr}

\author{F. Balacheff}

\address{F. Balacheff, Laboratoire Paul Painlev\'e, Bat. M2, Universit\'e des Sciences et Technologies, 59 655 Villeneuve d'Ascq, France.}

\email{florent.balacheff@math.univ-lille1.fr}

\author{K. Tzanev}

\address{K. Tzanev, Laboratoire Paul Painlev\'e, Bat. M2, Universit\'e des Sciences et Technologies, 59 655 Villeneuve d'Ascq, France.}

\email{kroum.tzanev@math.univ-lille1.fr}



\thanks{This work was partially supported by the grant ANR12-BS01-0009 FINSLER.}

\begin{document}%

\begin{abstract}
    We explore a natural generalization of systolic geometry to Finsler metrics and optical hypersurfaces with special emphasis on its relation to the Mahler conjecture and the geometry of numbers. In particular, we show that if an optical hypersurface of contact type in the cotangent bundle of the $2$-dimensional torus encloses a volume $V$, then it carries a periodic characteristic whose action is at most $\sqrt{V/3}$. This result is deduced from an interesting dual version of Minkowski's lattice-point theorem: if the origin is the unique integer point in the interior of a planar convex body, the area of its dual body is at least 3/2.
\end{abstract}

\maketitle

\begin{flushright}
    \begin{small}
        \parbox{4in} {\textit{Never consider a convex body without considering its dual at the same time.}} \\
        --- I.M. Gelfand
    \end{small}
\end{flushright}

\section{Introduction}

Minkowski's first theorem in the geometry of numbers states that \assertion{if the volume of a $0$-symmetric convex body in $\R^n$ is at least $2^n$, the body contains a non-zero integer point}. On the other hand, it is easy to find asymmetric convex bodies of arbitrary large volume that contain the origin and no other integer point. It is tempting to say something about the geometry of such bodies. For example, it is known that they must be flat in some lattice direction (see~\cite{Kannan-Lovasz:1988} and~\cite{Banaszczyk-etal:1999}). In this paper we show that the volume of their duals cannot be arbitrarily small. In fact, the interplay between contact and systolic geometry studied in~\cite{Alvarez-Balacheff:2011} suggests the following sharp inequality:

\begin{conjecture}{I}\label{intro:conj:I}
    If the interior of a convex body in $\R^n$ contains no integer point other than the origin, then the volume of its dual body is at least $(n+1)/n!$. Moreover, equality holds if and only if the convex body is a simplex such that the integer points on its boundary are precisely its vertices.
\end{conjecture}

Another formulation of the conjecture that seems more elementary is as follows: \assertion{if every integer hyperplane $m_1 x_1 + \cdots + m_n x_n = 1$---where the $m_i$ are integers not all equal to zero---intersects a convex body $K \subset \R^n$, then the volume of $K$ is at least $(n+1)/n!$}.

We prove both the two-dimensional case of the conjecture and its asymptotic version:

\begin{maintheorem}{I}\label{intro:thm:I}
    The area of a convex body in the plane that intersects every integer line $mx + ny = 1$ is at least $3/2$. Moreover, equality holds only for the triangle with vertices $(1,0)$, $(0,1)$, $(-1,-1)$ and its images under $GL(2,\Z)$.
\end{maintheorem}

\medskip
\begin{center}
    \includegraphics{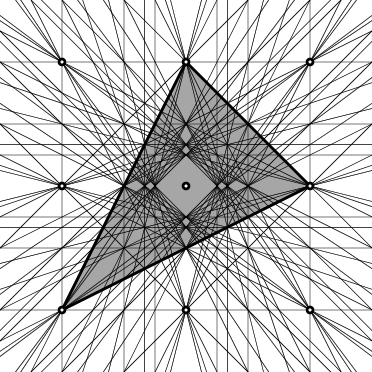}
    \captionof{figure}{Integer lines $mx + ny = 1$ with $-4 \leq m, n \leq 4$ and a convex body of minimal area that intersects every integer line.}
\end{center}
\medskip

\begin{maintheorem}{II}\label{intro:thm:II}
    The volume of a convex body in $\R^n$ that intersects every integer hyperplane $m_1 x_1 + \cdots + m_n x_n = 1$ is at least $(\pi/4)^n/n!.$ In particular the volume of such a convex body is at least
    $$
        c^n \, \frac{n+1}{n!}
    $$
    for any $c<\frac{\pi}{4}\simeq 0,785...$ provided that $n$ is large enough.
\end{maintheorem}

Theorem~\ref{intro:thm:II} is shown to be equivalent to the Bourgain-Milman theorem (\cite{Bourgain-Milman:1987} and~\cite{Kuperberg:2008}) in the sense that they imply each other in a relatively simple way. Despite this rough equivalence between Conjecture~\ref{intro:conj:I} and the Mahler conjecture, we believe the former to be somewhat simpler. For example, while it has long been an open problem to show that the minima of the volume product are polytopes (see~\cite{Reisner-Schutt-Werner:2012} for some recent progress on this front), we are able to prove the following

\begin{maintheorem}{III}\label{intro:thm:III}
    Among the convex bodies in $\R^n$ that intersect every integer hyperplane there exist bodies of minimal volume. These minimal bodies are polytopes and the number of their vertices is bounded by a quantity $c_n$ that depends only on the dimension $n$.
\end{maintheorem}

Our interest in these results stems from their relation to Hamiltonian dynamics. Indeed, Theorem~\ref{intro:thm:I} generalizes to a (sharp) isosystolic inequality for \term{optical hypersurfaces of contact type} on the cotangent bundle of the $2$-torus. These hypersurfaces, which we shall henceforth call optical hypersurfaces for brevity's sake, have a simple definition.

\begin{definition*}
    A smooth hypersurface in the cotangent bundle of a manifold $M$ is said to be \term{optical} if its intersection with each cotangent space is a quadratically convex hypersurface enclosing the origin. If all these convex hypersurfaces are $0$-symmetric, we shall say that the optical hypersurface is \term{reversible}.
\end{definition*}

\begin{maintheorem}{IV}\label{intro:thm:IV}
    If an optical hypersurface in the cotangent bundle of the $2$-torus bounds a volume $V$, it carries a periodic characteristic whose action is at most $\sqrt{V/3}$.
\end{maintheorem}

We advise the reader that in this theorem, as in the rest of the paper, volumes on the cotangent bundle of an $n$-dimensional manifold are measured with the $n$-th power of the symplectic form $\omega^n$. Note this is $n!$ times the usual definition of the symplectic or Liouville volume. Theorem~\ref{intro:thm:IV} can be restated in terms of Finsler geometry: \assertion{if the Holmes-Thompson volume of a Finsler $2$-torus is $3/2\pi$, then it carries a periodic geodesic of length at most 1}. This is the Finsler generalization of Loewner's isosystolic inequality. For \term{reversible} Finsler metrics the constant $3/2\pi$ can be improved to $2/\pi$, a result due to S.~Sabourau~(\cite{Sabourau:2010}) and one of the starting points of our investigation.

The reason for choosing the language of optical hypersurfaces over that of Finsler metrics (at least in this introduction) is our belief that some of these results will generalize to the contact-geometric setting for systolic geometry proposed in~\cite{Alvarez-Balacheff:2011}. Indeed, one of our aims in this paper is to show that optical hypersurfaces are a safe and interesting testing ground for this further contact-geometric generalization. On the one hand, convex-geometric inequalities allow us to prove the following extensions of key results of Gromov~\cite{Gromov:1983} and Croke~\cite{Croke:1988} to optical hypersurfaces:

\begin{maintheorem}{V}\label{intro:thm:V}
    Given an essential $n$-dimensional manifold $M$, there exists a constant $C > 0$ such that every optical hypersurface $\Sigma \subset T^* M$ bounding a volume $V$ carries a closed characteristic whose action is less than $C \sqrt[n]V$.
\end{maintheorem}

\begin{maintheorem}{VI}\label{intro:thm:VI}
    There exists constant $C > 0$ such that every optical hypersurface $\Sigma \subset T^* S^2$ bounding a volume $V$ carries a closed characteristic whose action is less than $C \sqrt{V}$.
    \footnote{
        The explicit value for $C$ coming out of our proof is huge. Since this paper appeared in the arXiv, the constant has been improved by Y.~Liokumovich in \cite{Liokumovich:2014} to $160$, and then by F.~Balacheff in \cite{Balacheff:2015} to $32$. Both improvements are based on the arguments presented in subsection \ref{coarse-inequalities-subsection} (more specifically Theorems~\ref{inequality-reversible} and~\ref{inequality-non-reversible}).
    }
\end{maintheorem}

On the other hand, the new setting strengthens the deep ties between systolic geometry, the geometry of numbers, and convex geometry. Conjecture~\ref{intro:conj:I} and the theorems around it illustrate this. Another illustration centers around the following straightforward generalization of the conjectured systolic optimality of the canonical metric in real projective spaces:

\begin{conjecture}{II}\label{intro:conj:II}
    If a reversible optical hypersurface in the cotangent bundle of real projective $n$-space encloses a volume $V$, it carries a periodic characteristic whose action is at most $\sqrt[n]{V}/2$. Moreover, this short characteristic can be chosen so that its projection onto the base manifold is a non-contractible loop.
\end{conjecture}

In \S\ref{Mahler-conjecture-subsection} we discuss this conjecture and show that Mahler's conjectured lower bound for the volume product of centrally symmetric convex bodies follows easily from it.

\begin{maintheorem}{VII}\label{intro:thm:VII}
    The conjectured sharp isosystolic inequality for reversible optical hypersurfaces in $T^*\RP^n$ implies the Mahler conjecture for centrally symmetric convex bodies in $\R^n$.
\end{maintheorem}

It is notable that, in conjunction with S.~Ivanov's proof of Conjecture~\ref{intro:conj:II} for $n = 2$ (see~\cite{Ivanov:2002} and~\cite{Ivanov:2011}), Theorem~\ref{intro:thm:VII} yields a new proof of Mahler's inequality for centrally symmetric convex bodies on the plane: \assertion{the product of the areas of a $0$-symmetric convex body on the plane and its dual body is greater than or equal to $8$}.

Our impression is that there are many other as-yet-undiscovered connections between convex geometry, the geometry of numbers, symplectic geometry, and isosystolic inequalities. In consequence, we have made an effort to make a large part of the paper accessible to readers from these different backgrounds.

\subsubsection{Plan of the paper}
Section~\ref{prelims} is a short review of important concepts and inequalities in convex geometry. The expert should just skim through it in order to become acquainted with some of the notation and terminology used in the rest of the paper. Our results on the geometry of numbers as well as a thorough comparison with previous results will be found in Section~\ref{geometry-of-numbers}. More precisely, the reader will find the proof of Theorem~\ref{thm:II} in \S\ref{Q-star}, the proof of Theorem~\ref{thm:III} in \S\ref{minima-subsection}, and two proofs of Theorem~\ref{thm:I} in \S\ref{quick-subsection} and \S\ref{algorithm-subsection}. The results on the systolic geometry of Finsler metrics and optical hypersurfaces will be found in Section~\ref{systolic}. The reader will find the proof of Theorem~\ref{thm:IV} in \S\ref{Loewner-subsection} and the proofs of Theorems~V and~VI in \S\ref{coarse-inequalities-subsection}. Theorem~\ref{thm:VII} is the subject of \S\ref{Mahler-conjecture-subsection}. The underlying connection between the results of this paper and the contact-geometric approach to systolic geometry in~\cite{Alvarez-Balacheff:2011} is elucidated in \S\ref{Zoll-tori-subsection}.

\newpage

\begin{center}
    \textbf{Acknowledgments}
\end{center}
\smallskip

The authors heartily thank Marie-Claude Arnaud, Ivan Babenko, Franck Barthe, Dario Cordero-Erausquin, Yaron Ostrover, Hugo Parlier, Gabriel Rivi\`ere, St\'ephane Sabourau, Tony Thompson, and Constantin Vernicos for some interesting conversations. We are also grateful for some very helpful MathOverflow exchanges with Noam Elkies, Sergei Ivanov, and Yoav Kallus. Special thanks go to E.~Makai~Jr. for stimulating exchanges and for making his work available to us. The authors are also grateful to the anonymous referee for the useful comments. Lastly, the epigraph to this paper is a piece of advice I.M.~Gelfand once gave J.-C. \'Alvarez Paiva: ``Never consider a convex body without considering its dual at the same time. If you consider both together, you have symplectic geometry and the uncertainty principle.'' The present work applies this advice to the geometry of numbers.

\section{A compendium of convex geometry}\label{prelims}

In this section we quickly review some basic notions and inequalities in convex geometry. For proofs and details we refer the reader to the books of A.C.~Thompson and P.~Gruber (\cite{Thompson:1996} and~\cite{Gruber}).

\subsection{Asymmetric norms}

Let $V$ be a finite-dimensional vector space over the reals and let $K \subset V$ be a convex body containing the origin as an interior point. The \term{gauge function} or (asymmetric) \term{norm } associated to $K$ is the function $\| \cdot \|_K : V \rightarrow [0, \infty)$ defined by
$$
    \|v\|_K := \inf \{t > 0 : v \in tK \} .
$$
As is well known, $\| \cdot \|_K$ satisfies the following properties:
\begin{enumerate}
    \item $\|v\| _K \geq 0$ with equality if and only if $v$ is the zero vector;
    \item $\|\lambda v\|_K = \lambda \|v\|_K$ for every positive real number $\lambda$;
    \item $\|v + w\|_K \leq \|v\| _K+ \|w\|_K$.
\end{enumerate}
Conversely, if $\| \cdot \|$ is a function satisfying these three properties, its \term{unit ball}
$$
    B := \{v \in V : \|v\| \leq 1\}
$$
is a convex body containing the origin as an interior point and $\| \cdot \| = \| \cdot \|_B$. Note that if $T : V \rightarrow V$ is an invertible linear transformation, then
$$
    \| T^{-1} v \|_B = \|v\|_{T(B)} \text{\quad for all } v \in V .
$$

Two special classes of norms that will be useful later in the paper are \term{reversible norms} (i.e., $\|-v\| = \|v\|$ for all $v \in V$) and \term{Minkowski norms}: outside the origin the function $\| \cdot \|^2$ is smooth and its Hessian is positive-definite. The unit sphere of a Minkowski norm is a \term{quadratically convex hypersurface}: the osculating quadric at each of its points is an ellipsoid.

\subsection{The space of convex bodies}

The space $\mathcal{K}_0(V)$ of all convex bodies in $V$ that contain the origin as an interior point will be topologized as a subset of the space $\mathcal{C}(V)$ of all compact subsets of $V$, which in turn is provided with the topology induced from the Hausdorff metric (see \cite[p.~65]{Thompson:1996}). In its modern formulation, Blaschke's selection theorem states that if $C \subset V$ is compact, then the set of all closed subsets of $C$ is itself a compact subset of $\mathcal{C}(V)$. This implies the following useful result:

\begin{proposition}\label{compactness}
    A closed subset $X \subset \mathcal{K}_0(V)$ is compact if and only if there exist two nested convex bodies $K_1 \subset K_2$ belonging to $\mathcal{K}_0(V)$ such that
    $$
        X \subset \{K \in \mathcal{K}_0(V): K_1 \subset K \subset K_2 \} .
    $$
\end{proposition}

The natural action of $GL(V)$, the group of invertible linear transformations from $V$ to itself, on the space $\mathcal{K}_0(V)$ is continuous and preserves the subspace $\mathcal{K}_0^s(V) \subset \mathcal{K}_0(V)$ of $0$-symmetric convex bodies. The main difference between the study of $0$-symmetric convex bodies and the study of convex bodies that contain the origin as an interior point is that the quotient space $\mathcal{K}_0^s(V)/GL(V)$---known as the \term{Banach-Mazur compactum}---is compact, while $\mathcal{K}_0(V)/GL(V)$ is not. In practical terms, this means that while every continuous, linear-invariant functional defined on $\mathcal{K}_0^s(V)$ is bounded and attains its extremal values, this will not be the case for $\mathcal{K}_0(V)$.

For continuous, \term{affine-invariant} functionals we have the following result of Macbeath~\cite{Macbeath:1951}:

\begin{theorem}[Macbeath]\label{Macbeath}
    The space of affine equivalence classes of convex bodies in a finite-dimensional (real) vector space is compact. In particular, every continuous affine invariant of convex bodies is bounded and attains its extremal values.
\end{theorem}

\subsection{Duality}

Given a norm on a finite-dimensional vector space $V$, we define its dual norm on $V^*$ by
$$
    \|\xi\|^* = \sup\{\xi \cdot v : \|v\| \leq 1 \}.
$$
If $B$ is the unit ball of $(V,\|\cdot\|)$, the unit ball of the dual normed space $(V^*,\|\cdot\|^*)$ will be denoted by $B^*$ and called the \term{dual body} of $B$. The following statements are standard interpretations of the dual body and its boundary:
\begin{itemize}
    \item $B^* \! \setminus \! 0$ is the set of hyperplanes that do not intersect the interior of $B$: $\xi \in B^* \! \setminus \! 0$ if and only if the hyperplane $\{v \in V: \xi \cdot v = 1\}$ does not intersect the interior of $B$.
    \item The dual unit sphere is the set of hyperplanes supporting $B$.
\end{itemize}

The duality operation on convex bodies satisfies some important properties, which we summarize in the following

\begin{proposition}
    The map $B \mapsto B^*$ is a continuous map from $\mathcal{K}_0(V)$ to $\mathcal{K}_0(V^*)$ and satisfies the following properties:
    \begin{enumerate}
        \item $(B^*)^* = B$;
        \item if $B \subset K$, then $K^* \subset B^*$;
        \item if $T: V \rightarrow V$ is an invertible linear map, then $T^{-1}(B)^* = T^*(B^*)$;
        \item $B$ is a polytope if and only if $B^*$ is a polytope;
        \item the hypersurface $\partial B$ is quadratically convex if and only if $\partial B^*$ is quadratically convex.
    \end{enumerate}
\end{proposition}

From property~(3) it follows that $B$ is symmetric about the origin if and only if $B^*$ is also symmetric and that dilating $B$ by a factor $\lambda > 0$ dilates its dual body by a factor $1/\lambda$. Property~(5) is equivalent to the statement: \assertion{the dual of a Minkowski norm is also a Minkowski norm}.

\subsection{Linear and affine invariants of convex bodies}

Fix a translation-invariant volume density on $V$ (i.e., a reversible norm on the one-dimensional space $\bigwedge^n V$) and consider its dual density on $V^*$ defined by the relation
$$
    |v_1^* \wedge \cdots \wedge v_n^*| = |v_1 \wedge \cdots \wedge v_n|^{-1}
$$
whenever $v_1, \ldots, v_n$ is a basis of $V$ and $v_1^*, \ldots, v_n^*$ is the dual basis in $V^*$. Note that if $U$ is any Borel subset of $V$ and $|v_1 \wedge \cdots \wedge v_n| = 1$, the integral of the volume density over $U$---to be denoted by $|U|$---is nothing more than the ratio of the measure of $U$ and the volume of the parallelotope formed by the vectors $v_1, \ldots, v_n$ when both quantities are computed with a fixed, but otherwise arbitrary, Lebesgue measure on $V$.

The two basic linear invariants of a convex body $K \in \mathcal{K}_0(V)$ are
$$
    |K||K^*| \text{\quad and\quad} |K - K|/|K| .
$$
Note that $|K||K^*|$ is $1/n!$ times the volume of $K \times K^* \subset V \times V^*$ measured with the $n$-th exterior power of the standard symplectic form on $V \times V^*$, while $|K - K|/|K|$ is the relative volume of the \term{difference body}
$$
    K - K := \{ x - y \in V : x, y \in K\}
$$
and the body $K$.

The Brunn-Minkowski and Rogers-Shephard inequalities (see~\cite{Rogers-Shephard:1957}) determine the minimum and maximum values of $|K - K|/|K|$:

\begin{theorem}\label{Rogers-Shephard}
    For any convex body $K \subset V$, we have that
    $$
        2^n \leq \frac{|K - K|}{|K|} \leq \frac{(2n)!}{(n!)^2} .
    $$
    Moreover, the left-hand side becomes an equality if and only if $K$ is symmetric about some point, whereas the right-hand side becomes an equality if and only if $K$ is a simplex.
\end{theorem}

In the following theorem, the inequality on the left is G.~Kuperberg's remarkable sharpening of the Bourgain-Milman theorem (see Kuperberg~\cite{Kuperberg:2008}), whereas on the right is the well-known Blaschke-Santal\'o inequality.

\begin{theorem}\label{Santalo-Kuperberg}
    For any convex body $K \subset V$, we have that
    $$
        \frac{\pi^n}{n!} < |K-K||(K-K)^*| \leq \varepsilon_n^2 ,
    $$
    where $\varepsilon_n$ is the volume of the $n$-dimensional Euclidean unit ball. Moreover, equality on the right holds if and only if $K-K$ is an ellipsoid.
\end{theorem}

The sharp lower bounds for $|K||K^*|$ and $ |K-K||(K-K)^*|$ are still conjectural in dimension greater than two:

\begin{theorem}[Mahler~\cite{Mahler-theorem:1939}] \label{Mahler-theorem}
    If $V$ is two-dimensional and $K \subset V$ is a convex body with the origin in its interior,
    $$
        \frac{27}{4} \leq |K||K^*| \text{\quad and\quad} 8 \leq |K-K||(K-K)^*| .
    $$
    The first inequality becomes an equality if and only if $K$ is a triangle with the origin as barycenter, whereas the second inequality becomes an equality if and only if $K -K$ is a parallelogram.
\end{theorem}

\begin{mahlerconjecture}[\cite{Mahler-conjecture:1939}]\label{Mahler-conjecture}
    For every $K \in \mathcal{K}_0(V)$
    $$
        \frac{(n+1)^{n + 1}}{(n!)^2} \leq |K| |K^*| \text{\quad and\quad} \frac{4^n}{n!} \leq |K-K||(K-K)^*| .
    $$
\end{mahlerconjecture}

The reader will have already noticed that given a continuous, linear-invariant functional $\mathcal{F}(K)$, we may construct the symmetrized invariant $\mathcal{F}(K-K)$. This is an idempotent operation on the space of continuous functions on $\mathcal{K}_0(V)/GL(V)$. It is noteworthy that $\mathcal{F}(K-K)$ is naturally a continuous affine invariant on the space of all convex bodies in $V$. By Macbeath's compactness theorem (Theorem~\ref{Macbeath}), this functional is bounded and attains is extremal values.

Besides this idempotent operation on the space of continuous, linear-invariant functionals on $\mathcal{K}_0(V)$, we also have a natural involution that is suggested by the theory of volumes on normed spaces:

\begin{definition}
    If $\mathcal{F} : \mathcal{K}_0(V) \rightarrow (0,\infty)$ is a linear-invariant functional, we define the \term{dual functional} $\mathcal{F}^*$ by the relation
    $$
        \mathcal{F}^*(K) \mathcal{F}(K^*) = |K||K^*| .
    $$
\end{definition}

This definition, while somewhat intriguing, is not as arbitrary as it may seem and is closely related to the notion of dual functor on the category of normed spaces (see~\cite[p.~11]{Alvarez-Thompson:2004} and~\cite{Alvarez-Thompson:Monthly}).

From a positive linear invariant $\mathcal{F}(K)$ we may derive, besides the linear invariant $\mathcal{F}^*(K)$, the affine invariants
\begin{itemize}
    \item $\mathcal{F}(K-K)$,
    \item $\mathcal{F}^*(K-K)$,
    \item $\mathcal{F}(K-K)(|K|/|K-K|)$,
    \item $\mathcal{F}^*(K-K)(|K|/|K-K|)$.
\end{itemize}
Of course, other invariants can be obtained by this formal game involving symmetrization and duality. We would just like to call the reader's attention to a simple principle that will allow us to organize many of the classic and not-so-classic results in the geometry of numbers, as well as to provide context to the results presented in this paper: through Theorems~\ref{Rogers-Shephard} and~\ref{Santalo-Kuperberg}, any bounds for the invariant $\mathcal{F}$ \term{automatically} translate into bounds for the derived invariants.

\section{Results in the geometry of numbers}\label{geometry-of-numbers}

In this section, we are interested in studying the action of the group $GL(n,\Z)$ of unimodular transformations on the space $\mathcal{K}_0(\R^n)$ of convex bodies in $\R^n$ that contain the origin as an interior point. This is a classical subject which is part of \term{the geometry of numbers} and is home to some of the most beautiful theorems in geometry. Among these we find Minkowski's lattice-point theorem and the Minkowski-Hlawka theorem:

\begin{theorem}[Minkowski]\label{Minkowski}
    If the volume of a $0$-symmetric convex body in $\R^n$ is at least $2^n$, the body contains a non-zero integer point.
\end{theorem}

\begin{theorem}[Minkowski and Hlawka]\label{Minkowski-Hlawka}
    Let $K \subset \R^n$ be a convex body containing the origin as an interior point. If $|K| < 1$, there exists a unimodular transformation $T$ for which the body $T(K)$ does not contain any non-zero integer point. Moreover, if $K$ is $0$-symmetric, we obtain the same conclusion with the weakened hypothesis $|K| < 2$.
\end{theorem}

For proofs of these results---and a more general version of the Minkowski-Hlawka theorem---we refer the reader to the books~\cite{Cassels:1971, Lekkerkerker:1969} and~\cite{Gruber}, which will also be our basic references for the rest of this section.

\subsection{The invariant \texorpdfstring{$\mathcal{Q}$}{Q} }

It has proved useful to inscribe the two preceding theorems in the framework of the study of linear and affine invariants of convex bodies. This is done by defining the \term{critical determinant} of a convex body $K \in \mathcal{K}_0(\R^n)$ as the minimum of the determinants of all lattices that intersect its interior only at the origin. The critical determinant, denoted by $\Delta$, is not itself a linear invariant, but it may be used to define the invariant
$$
    \mathcal{Q}(K) := \frac{|K|}{\Delta(K)}
$$
(the invariance will follow from Proposition~\ref{critical-determinant}). The following translation of Theorems~\ref{Minkowski} and~\ref{Minkowski-Hlawka} follows easily from the linear invariance of $\mathcal{Q}$ and the fact that every lattice in $\R^n$ is linearly equivalent to $\Z^n$.

\begin{theorem}\label{Minkowski-Hlawka-inequality}
    For every $K \in \mathcal{K}_0(\R^n)$,
    $$
        1 \leq \mathcal{Q}(K) \text{\quad and\quad} 2 \leq \mathcal{Q}(K-K) \leq 2^n .
    $$
\end{theorem}

All the properties of the critical determinant that we will need are contained in the following

\begin{proposition}~\label{critical-determinant}
    The functional $\Delta : \mathcal{K}_0(\R^n) \rightarrow (0, \infty)$ satisfies the following properties:
    \begin{enumerate}
        \item if $K_1 \subset K_2$, then $\Delta(K_1) \leq \Delta(K_2)$;
        \item $\Delta(T(K)) = |\det(T)| \Delta(K)$ for any invertible linear map $T$;
        \item $2^n \Delta((K - K)^*) \leq \Delta(K^*)$.
    \end{enumerate}
\end{proposition}

\begin{proof}
    The first two properties follow easily from the definition of the critical determinant. Together, they imply that $\Delta$ is a continuous functional on $\mathcal{K}_0(\R^n)$ (see~\cite[\S~25.3]{Lekkerkerker:1969}).

    In order to verify the third property, we first rewrite it as
    $$
        \Delta \left(\left( \frac{K-K}{2}\right)^*\right) \leq \Delta(K^*)
    $$
    and remark that for $\xi \in \R^{n*}$
    $$
        \|\xi\|_{\left( \frac{K-K}{2}\right)^*} = \frac{1}{2} \left(\|\xi\|_{K^*} + \|-\xi\|_{K^*} \right) .
    $$
    This implies that the set of lattices
    $$
        \{ \Lambda^* \subset \R^{n*} : 1 \leq \|\xi\|_{\left( \frac{K-K}{2}\right)^*} \text{ for } \xi \in \Lambda^* \! \setminus \! \{ 0 \} \}
    $$
    contains $\{ \Lambda^* \subset \R^{n*} : 1 \leq \|\xi\|_{K^*} \text{ for } \xi \in \Lambda^* \! \setminus \! \{ 0 \} \} $ and, therefore, the desired inequality follows from the definition of the critical determinant.
\end{proof}

Almost every quantity studied in the geometry of numbers can be written in terms of the critical determinant. For example, another classic result of Minkowski (see~\cite{Minkowski:1904} and~\cite[p.~69]{Rogers:1964}) states that the density of the tightest lattice packing of a convex body $K \subset \R^n$ is given by
$$
    \delta(K) := |K|/\Delta(K-K) = \mathcal{Q}(K-K)\frac{|K|}{|K-K|} .
$$
Notice that the inequalities $2(n!)^2/(2n)! \leq \delta(K) \leq 1$ follow at once from Theorems~\ref{Minkowski-Hlawka} and~\ref{Rogers-Shephard}. Indeed, applying the simple principle presented at the end of Section~\ref{prelims}, we see that the packing density $\delta$ is but one of a host of linear and affine invariants derived from $\mathcal{Q}$ and which \term{automatically} satisfy certain inequalities implied by Theorems~\ref{Minkowski-Hlawka-inequality}, \ref{Rogers-Shephard}, and~\ref{Santalo-Kuperberg}. What is remarkable about the invariant $\mathcal{Q}$ is that most, if not all, of its derived invariants have interesting and non-trivial geometric interpretations.

\subsection{The invariant \texorpdfstring{$\mathcal{Q}^*$}{Q*} }\label{Q-star}

In this paper, our main interest is the dual invariant $\mathcal{Q}^*(K) = |K| \Delta(K^*)$. Given that every lattice in $\R^n$ is linearly equivalent to $\Z^n$, the linear invariance of $\mathcal{Q}^*$ allows us to rewrite
\begin{conjecture}{I}\label{conj:I}
    If the interior of a convex body in $\R^n$ contains no integer point other than the origin, then the volume of its dual body is at least $(n+1)/n!$. Equivalently,
    $$
        \frac{n+1}{n!} \leq \mathcal{Q}^*(K) \text{\quad for all } K \in \mathcal{K}_0(\R^n).
    $$
    Moreover, equality holds if and only if the convex body is a simplex such that the integer points on its boundary are precisely its vertices.
\end{conjecture}

Although it seems that the invariant $\mathcal{Q}^*$ is considered here for the first time, the invariant $\mathcal{Q}^*(K-K)$ was already considered by Mahler in~\cite{Mahler:1974} and~\cite{Mahler:1976}, and was later studied in Makai~\cite{Makai:1978}. Apart from the application of G.~Kuperberg's lower bound for $|K-K||(K-K)^*|$ and a somewhat different geometric interpretation, the following result is due to Mahler (\cite{Mahler:1974, Mahler:1976}).

\begin{proposition}\label{Mahler-dual}
    If $K \subset \R^n$ is a convex body,
    $$
        |(K-K)^*||K-K| \leq 2^n \mathcal{Q}^*(K-K) .
    $$
    Consequently, the volume of a $0$-symmetric convex body $B \subset \R^n$ that intersects every integer hyperplane is greater than $(\pi/2)^n/ n!$. Moreover, if $n = 2$, then $2 \leq |B|$ with equality if and only if $B$ is a parallelogram.
\end{proposition}

\begin{proof}
    By Minkowski's lattice-point theorem $\mathcal{Q}((K-K)^*) \leq 2^n$ and, therefore,
    $$
        |(K-K)^*||K-K| = \mathcal{Q}((K-K)^*) \mathcal{Q}^*(K-K) \leq 2^n \mathcal{Q}^*(K-K) .
    $$
    The second part of the result follows from Theorems~\ref{Santalo-Kuperberg} and~\ref{Mahler-theorem}.
\end{proof}

In the same way that the packing density $\delta$ was defined in terms of $\mathcal{Q}$, we may define the invariant
$$
    \rho(K) := |K| \Delta((K-K)^*) = \mathcal{Q}^*(K-K) \frac{|K|}{|K-K|} .
$$
The main result of Makai~\cite{Makai:1978} identifies this invariant as the density of the thinnest lattice of translates of $K$ that intersects every affine hyperplane. The relation between $\mathcal{Q}^*(K)$ and this invariant is given by the following

\begin{proposition}\label{inequalities}
    For every $K \in \mathcal{K}_0(\R^n)$,
    $$
        2^n \rho(K) \leq \mathcal{Q}^*(K) < |K||K^*|.
    $$
\end{proposition}

\begin{proof}
    Dividing all terms in the inequality by $|K|$ reduces to proof to showing that
    $$
        2^n \Delta((K-K)^*) \leq \Delta(K^*) \leq |K^*| .
    $$
    The inequality on the left follows immediately from Proposition~\ref{critical-determinant}, while the inequality on the right follows immediately from the Minkowski-Hlawka inequality.
\end{proof}

We are now ready to prove

\begin{maintheorem}{II}\label{thm:II}
    The volume of a convex body in $\R^n$ that intersects every integer hyperplane $m_1 x_1 + \cdots + m_n x_n = 1$ is at least $(\pi/4)^n /n!$.
\end{maintheorem}

The first step will be to show that an \term{unavoidable convex body}---a convex body that intersects every integer hyperplane---contains the origin as an interior point. This will also allow us to characterize unavoidable convex bodies as those whose dual bodies do not contain any non-zero integer point in the interior.

\begin{lemma}\label{origin-inside}
    A convex body $K$ contains the origin as an interior point if and only if there is at most a finite number of integer hyperplanes that do not intersect $K$. In particular, every unavoidable convex body contains the origin in its interior.
\end{lemma}

\begin{proof}
    Since the distance between the hyperplane $m_1x_1 + \cdots + m_nx_n = 1$ and the origin is $1/\sqrt{m_1^2 + \cdots + m_n^2}$, there are only a finite number of integer hyperplanes that do not intersect a given neighborhood of the origin.

    We now show that if the origin is not in the interior of a convex body $K$, it must fail to intersect an infinite number of integer hyperplanes.

    If we assume the convex body $K$ does not contain the origin in its interior, then $K$ is contained in a closed half-space $\xi \cdot x \leq 0$, where $\xi$ is a nonzero element of $\R^{n*}$. Let $v$ be any vector satisfying $\xi \cdot v < 0$ and consider the projection $\Pi : \R^n \rightarrow \R^n$ onto the hyperplane $\xi = 0$ along the $1$-dimensional subspace spanned by $v$. Let us enclose $\Pi(K)$ in an $(n-1)$-dimensional euclidean disc $D$ lying on the hyperplane $\xi = 0$ and centered at the origin. The dual body of $D$ in $\R^{n*}$ is a $0$-symmetric cylinder (of infinite volume) and hence, by van~der~Corput's refinement of Minkoswki's lattice-point theorem (see~\cite[\S7.2]{Lekkerkerker:1969}), it contains an infinite number of non-zero integer points $\pm \eta_i$ $(i \in \N)$. It follows that the integer hyperplanes $\eta_i = \pm 1$ $(i \in \N)$ do not intersect $D$ nor $\Pi(K)$. In other words, the projection of $K$ is contained in every slab $-1 < \eta_i < 1$ $(i \in \N)$.

    Assuming, without loss of generality, that $\eta_i \cdot v \geq 0$ for all $i \in \N$, we remark that the integer hyperplanes $\eta_i = -1$ do not intersect $K$.
\end{proof}

\begin{proof}[Proof of Theorem~\ref{thm:II}]
    If $K \subset \R^n$ is an unavoidable convex body, the preceding lemma tells us that $K \in \mathcal{K}_0(\R^n)$ and, therefore, lies in the domain of the invariant $\mathcal{Q}^*$. Thus Theorem~\ref{thm:II} amounts to prove that
    $$
        \left(\frac{\pi}{4}\right)^n \, \frac{1}{n!} \leq \mathcal{Q}^*(K).
    $$
    By Minkowski's lattice-point theorem, $|(K-K)^*| \leq 2^n \Delta((K-K)^*)$. Multiplying both sides by $|K|$ and applying Proposition~\ref{inequalities}, we obtain
    $$
        |K||(K-K)^*| \leq 2^n |K|\Delta((K-K)^*) \leq \mathcal{Q}^*(K) .
    $$
    Writing $|K||(K-K)^*|$ as $|K - K| |(K - K)^*| (|K|/|K-K|)$, Rogers-Shephard inequality in Theorem~ \ref{Rogers-Shephard} and Kuperberg's inequality in Theorem~\ref{Santalo-Kuperberg} imply that
    $$
        \frac{\pi^n n!}{(2n)!} \leq \mathcal{Q}^*(K) .
    $$
    The result now follows from the standard inequality
    $$
        \left(\frac{1}{4}\right)^n < \frac{(n!)^2}{(2n)!}
    $$
    that holds for every integer $n \geq 1$.
\end{proof}

We end the section by showing that Conjecture~\ref{conj:I} is roughly equivalent to the Mahler conjecture:

\begin{proposition}\label{ABT-implies-BM}
    If the constant $c > 0$ is such that
    $$
        c^n \, \frac{n+1}{n!} \leq \mathcal{Q}^*(K)
    $$
    for every convex body $K \in \mathcal{K}_0(\R^n)$, then
    $$
        \left(\frac{c}{e}\right)^n \frac{(n+1)^{n+1}}{(n!)^2} \leq |K||K^*|
    $$
    for every convex body $K \in \mathcal{K}_0(\R^n)$.
\end{proposition}

\begin{proof}
    Every convex body $K \in \mathcal{K}_0(\R^n)$ satisfies $\mathcal{Q}^*(K)<|K||K^*|$ according to Proposition~\ref{inequalities}. If $c^n \, \frac{n+1}{n!} \leq \mathcal{Q}^*(K)$ we get
    $$
        c^n \, \frac{n!}{(n+1)^{n}} \frac{(n+1)^{n+1}}{(n!)^2}< |K||K^*|.
    $$
    From the following refinement of Stirling's formula (see~\cite[p.~54]{Feller})
    $$
        \sqrt{2\pi n} \left(\frac{n}{e}\right)^n e^{1/(12n +1)} < n! ,
    $$
    we derive that
    $$
        \left(\frac{c}{e}\right)^n \frac{(n+1)^{n+1}}{(n!)^2} \leq |K||K^*|.
    $$
\end{proof}

\subsection{Existence and nature of the minima of \texorpdfstring{$\mathcal{Q}^*$}{Q*} }\label{minima-subsection}

Some of the functionals considered in this section such as $\mathcal{Q}(K-K)$, $\mathcal{Q}^*(K-K)$, and $\rho(K)$ are naturally continuous, affine-invariant functionals on the space $\mathcal{K}(\R^n)$ of all convex bodies in $\R^n$. This implies that they are bounded and attain their extremal values. Contrariwise, the continuous, linear-invariant functionals $|K||K^*|$, $\mathcal{Q}(K)$, and $\mathcal{Q}^*(K)$ defined on $\mathcal{K}_0(\R^n)$ are unbounded from above and the existence of minimal convex bodies must be determined. For $|K||K^*|$ this seems to be folklore, while it is not clear that there exist convex bodies that minimize $\mathcal{Q}(K)$. In what follows we show that there exist convex bodies in $\mathcal{K}_0(\R^n)$ that minimize the functional $\mathcal{Q}^*$ and that such bodies must be polytopes.

The main step in the proof of Theorem~\ref{thm:III} is the following

\begin{theorem}\label{compactness-theorem}
    The set of unavoidable convex bodies in $\R^n$ whose volume is bounded above by some constant $c > 0$ is compact modulo unimodular transformations. Equivalently, the sublevel sets of the functional $\mathcal{Q}^*$ are compact in $\mathcal{K}_0(\R^n)/GL(\R^n)$.
\end{theorem}

In turn, this result will follow from results of Minkowski and Mahler (see~\cite{Mahler:1938} and~\cite{Weyl:1940}) on the theory of reduction.

\begin{theorem}[Minkowski, Mahler]\label{reduction}
    Let $\| \cdot \| : \R^n \rightarrow [0, \infty)$ be a reversible norm and let $e_1,\ldots,e_n$ be the standard basis of $\R^n$. There exists a unimodular transformation $T$ such that
    $$
        a_i := \|T(e_i) \| = \min \{\| T(x)\|: x = (x_1,\ldots,x_n) \in \Z^n, \gcd(x_i,\ldots,x_n) = 1 \} .
    $$
    Moreover, if $B$ denotes the unit ball of the norm $\| \cdot \|$, the quantities $a_1, \ldots , a_n$ satisfy the inequality
    $$
        a_1 \cdots a_n |B| \leq 2^n \left(\frac{3}{2}\right)^{\frac{(n-1)(n-2)}{2}} =: \mu_n .
    $$
\end{theorem}

\begin{lemma}
    If $K \subset \R^n$ is any $0$-symmetric convex body that does not contain any non-zero integer point in its interior, there exists a unimodular transformation $T$ for which $T^{-1}(K)$ contains the ball with center at the origin and radius $|K|/\mu_n \sqrt{n}$.
\end{lemma}

\begin{proof}
    We apply Theorem~\ref{reduction} to the reversible norm with unit ball $K$ and remark that the definition of the quantities $a_i$ $(1 \leq i \leq n)$ implies that the convex body $a_n T^{-1}(K)$ contains the vectors $\pm e_1,\ldots, \pm e_n$. It follows that $a_n T^{-1}(K)$ also contains their convex hull, the cross-polytope $\ch\{\pm e_1,\ldots,\pm e_n\}$, and thus $T^{-1}(K)$ contains the ball with center at the origin and radius $1/a_n \sqrt{n}$.

    Since $T^{-1}(\mathring{K}) \cap \Z^n = \{0\}$, we have that $1 \leq a_i$ for $1 \leq i \leq n$ and, therefore,
    $$
        a_n |K| \leq a_1 \cdots a_n |K| \leq \mu_n \text{\quad or \quad} \frac{|K|}{\mu_n} \leq \frac{1}{a_n} .
    $$
    We conclude that $ T^{-1}(K)$ contains the ball with center at the origin and radius $|K|/\mu_n \sqrt{n}$.
\end{proof}

\begin{proof}[Proof of Theorem~\ref{compactness-theorem}]
    The first step in the two-step proof is to show that for every unavoidable convex body $K \subset \R^n$ there exists a unimodular transformation $T$ such that $T(K)$ is contained in a ball with center at the origin and radius $R_n |K|$, where $R_n > 0$ is a quantity that depends only on the dimension $n$. The second step is to show that the set of unavoidable convex bodies in $\R^n$ that are contained in a ball of fixed radius $R > 0$ is a compact subset of $\mathcal{K}_0(\R^n)$.

    \textit{First step.}
    Let us first assume that $K$ is $0$-symmetric. Applying the previous lemma to the body $K^*$ we obtain that there exists a unimodular transformation $T^*$ such that $T^{*-1}(K^*)$ contains the ball with center at the origin and radius $|K^*|/\mu_n \sqrt{n}$. By duality, the body
    $$
        \left(T^{*-1}(K^*)\right)^* = T(K)
    $$
    is contained in the ball with center at the origin and radius $\mu_n \sqrt{n}/|K^*|$. Using that $|K||K^*|$ is greater than some quantity $\gamma_n > 0$ that depends only on $n$ (we could use Theorem~\ref{Santalo-Kuperberg}, but any easier inequality would do), we conclude that $T(K)$ is contained in a ball centered at the origin and radius $R'_n |K| = \mu_n \sqrt{n}|K|/\gamma_n$.

    If the body $K$ is not $0$-symmetric, we apply the preceding argument to its difference body $K-K$. We have then that $K$ is contained in $K-K$ which in turn is contained in a ball centered at the origin and radius $\mu_n \sqrt{n}|K - K|/\gamma_n$. Since $|K-K| \leq (2n)! |K|/(n!)^2$, it follows from the Rogers-Shephard inequality that $K$ is contained in a ball centered at the origin and radius $R_n |K| = R'_n |K|(2n)!/(n!)^2$.

    \textit{Second step.}
    Since the topology of $\mathcal{K}_0(\R^n)$ is inherited from that of $\mathcal{C}(\R^n)$, the space of compact subsets of $\R^n$, it will be enough to prove that the set of unavoidable convex bodies in $\R^n$ that are contained in a ball of fixed radius $R > 0$ is compact in $\mathcal{C}(\R^n)$. This follows from the Blaschke selection theorem and the fact that intersecting all integer hyperplanes is a closed condition.
\end{proof}

\begin{maintheorem}{III}\label{thm:III}
    Among the convex bodies in $\R^n$ that intersect every integer hyperplane there exist bodies of minimal volume. Equivalently, the functional $\mathcal{Q}^*$ attains its minimum value on $\mathcal{K}_0(\R^n)$. Moreover, minimal bodies are polytopes and the number of their vertices is bounded by a quantity $c_n$ that depends only on the dimension $n$.
\end{maintheorem}

\begin{proof}
    The existence of minima follows from Theorem~\ref{compactness-theorem} and the continuity of $\mathcal{Q}^*$. To show that the minima are polytopes it suffices to show that \assertion{every unavoidable convex body contains an unavoidable polytope}.

    By Lemma~\ref{origin-inside}, in the interior of a unavoidable convex body $K$ we can find a non-degenerate simplex $S \in \mathcal{K}_0(\R^n)$. Notice that there is only a finite number of integer hyperplanes $H_1, \ldots, H_k$ not intersecting the interior of $S$. For each $i \in \{1,\ldots,k\}$, let $x_i$ be a point in the intersection $K \cap H_i$. The convex hull of the union of $S$ and $\{x_1,\ldots,x_k\}$ is an unavoidable polytope contained in $K$ and having no more than $n + 1 + k$ vertices.

    Since the volume of minimal unavoidable bodies is bounded above by $(n+1)/n!$, we can make use of the group of unimodular transformations---as in the first step of the proof of Theorem~\ref{compactness-theorem}---to restrict our attention to unavoidable polytopes contained in a ball of radius $R_n (n+1)/n!$. This is a compact subset of $\mathcal{K}_0(\R^n)$ and, by Proposition~\ref{compactness}, all unavoidable bodies in this set contain a common non-degenerate simplex in $\mathcal{K}_0(\R^n)$. It follows that there is a uniform bound on the number of vertices for all minimal unavoidable polytopes in $\R^n$.
\end{proof}

\subsection{A quick proof of Theorem~\ref{intro:thm:I}}\label{quick-subsection}

Let us remind the reader that Theorem~\ref{intro:thm:I} is the particular case of Conjecture~\ref{conj:I} in dimension $n=2$.
\begin{maintheorem}{I}\label{thm:I}
    The area of a convex body in the plane that intersects every integer line $mx + ny = 1$ is at least $3/2$. Equivalently,
    $$
        \frac{3}{2} \leq \mathcal{Q}^*(K) \text{\quad for all } K \in \mathcal{K}_0(\R^2).
    $$
    Moreover, equality holds only for the triangle with vertices $(1,0)$, $(0,1)$, $(-1,-1)$ and its images under $GL(2,\Z)$.
\end{maintheorem}
Our first proof of this theorem is a simple application of the following beautiful result of L.~Fejes-T\'oth and E.~Makai~Jr.~(\cite{Makai:1974}).

\begin{theorem}[L.~Fejes-T\'oth and E.~Makai~Jr.]
    If $K \subset \R^2$ is a convex body such that the set $K + \Z^2$ intersects every line in the plane, then the area of $K$ is at least $3/8$. Moreover equality holds if and only if, up to translations, $K$ is a triangle spanned by one vertex and the midpoints of the opposite sides of a basic parallelogram of the lattice $\Z^2$.
\end{theorem}

Indeed, this theorem is equivalent to the inequality $3/8 \leq \rho(K)$ and, therefore, Proposition~\ref{inequalities} implies that
$$
    \frac{3}{2} \leq 2^2 \rho(K) \leq \mathcal{Q}^*(K) \text{\quad for every } K \in \mathcal{K}_0(\R^2) .
$$

The authors hit upon this quick proof of Theorem~\ref{thm:I} after having found the algorithmic proof we will present in~\S\ref{algorithm-subsection}. However, we warn the reader that its ease is only apparent. To reconstruct a proof of Theorem~\ref{thm:I} with full details along these lines requires going through the results of three different publications: Makai~\cite{Makai:1978} for the proof that $\rho(K) = |K|\Delta((K-K)^*)$ is the density of the thinnest lattice of convex bodies that intersects every affine hyperplane; Fejes-T\'oth~and~Makai~\cite{Makai:1974} for the proof that, in the case of planar bodies, this density is at least $3/8$; and Behrend~\cite{Behrend:1937} for the reverse affine isodiametric inequality on which this last proof depends.

We mention in passing that in~\cite{Makai:1978} Makai conjectures that
$$
    \frac{n+1}{2^n n!} \leq \rho(K) \text{\quad for all convex bodies } K \subset \R^n .
$$
Thus, by Proposition~\ref{inequalities}, an affirmative solution to Makai's conjecture would immediately lead to an affirmative solution of Conjecture I. The geometrical fact behind the left-hand side inequality in Proposition~\ref{inequalities} simply states that if $K \subset \R^n$ is a convex body that intersects every integer hyperplane, then $\frac{1}{2}K + \Z^n$ intersects every hyperplane in $\R^n$.

\subsection{Algorithmic proof of Theorem~~\ref{thm:I}}\label{algorithm-subsection}

As we mentioned briefly in the introduction, the conjecture stating that the volume of an unavoidable convex body in $\R^n$ is at least $(n+1)/n!$ was suggested by the contact-geometric study of isosystolic inequalities in~\cite{Alvarez-Balacheff:2011}. In fact, this last paper also suggests a heuristic principle which, when translated to our present setting, implies that unavoidable bodies of minimal volume must be integer polytopes. An approach to the solution of the conjecture becomes apparent: (1) verify the conjecture for integer polytopes, (2) show that any unavoidable polytope $P \subset \R^n$ can be deformed to an unavoidable integer polytope in such a way that the volume does not increase along the deformation. Encouragingly, the first part of this approach is quite simple:

\begin{lemma}~\label{integer-polytope}
    The volume of an unavoidable integer polytope $P \subset \R^n$ is at least $(n+1)/n!$. Moreover, equality holds if and only if $P$ is a \term{basic simplex}: the image under an unimodular transformation of the simplex with vertices
    $$
        \big\{ (1,0 \ldots, 0), \ldots, (0,\ldots,0,1), (-1,-1,\ldots,-1) \big\}.
    $$
\end{lemma}

\begin{proof}
    We divide our proof into the three following claims.

    \begin{claim}
        An integer polytope is unavoidable if and only if it contains the origin as interior point.
    \end{claim}

    \begin{claimproof}
        By Lemma \ref{origin-inside}, any unavoidable convex body contains the origin in its interior. Conversely, let $P$ be an integer polytope containing the origin as interior point and consider any function $f(x_1,\ldots,x_n)=m_1x_1 + \cdots + m_nx_n$ where the $m_i$ are integers not all equal to zero. As the origin is an interior point of $P$ there must be at least one vertex $v \in P$ for which $f(v)>0$. Because $v \in \Z^n$, $f(v)$ is a non-zero integer. Thus $f(v)\geq 1$, and the integer hyperplane $m_1x_1 + \cdots + m_nx_n=1$ intersects $P$, proving that $P$ is indeed unavoidable.
    \end{claimproof}

    \begin{claim}
        The volume of an integer polytope $P \subset \R^n$ containing the origin as interior point is at least $(n+1)/n!$.
    \end{claim}

    \begin{claimproof}
        In each facet of the integer polytope $P$ we can find an $(n-1)$-dimensional simplex whose join with the origin is a non-degenerate lattice simplex. This construction can be performed at every facet of $P$ to obtain at least $n+1$ non-degenerate lattice simplices with disjoint interiors and whose union is contained in $P$. Since the volume of a non-degenerate lattice simplex is at least $1/n!$, the volume of $P$ is at least $(n+1)/n!$.
    \end{claimproof}

    \begin{claim}
        The volume of an integer polytope $P \subset \R^n$ containing the origin as interior point equals to $(n+1)/n!$ if and only if $P$ is a basic simplex.
    \end{claim}

    \begin{claimproof}
        If the volume of an integer polytope $P$ containing the origin as interior point equals $(n+1)/n!$, it must have exactly $n+1$ facets and must, therefore, be a lattice simplex containing the origin as interior point. Moreover, the simplices obtained as the convex hull of any one of its facets and the origin must have all volume $1/n!$. Thus the vertices of any of its facets form a basis of the integer lattice, which is equivalent for $P$ to be a basic simplex. Conversely, any basic simplex has volume $(n+1)/n!$.
    \end{claimproof}
\end{proof}

The part of the approach that calls for deforming the polytope is much more delicate. The first idea that comes to mind is to apply a steepest descent algorithm to the volume functional constrained to the space of unavoidable polytopes. For this idea to work, we would at the very least need to know that all local minima are integer polytopes. Nevertheless, in two-dimensions an optimization algorithm does give the desired answer. We describe its framework in very simple, abstract setting so that it may serve the reader as a blueprint for the more detailed constructions that follow.

\subsubsection{Abstract description of the algorithm}
Given a set $X$, a subset $Y \subset X$ and a function $f: X \rightarrow \R$, we are asked to determine the smallest value of $f(x)$ subject to the condition $x \in Y$. Evidently, establishing whether this smallest value exists and the set of points where it is attained is part of the problem. We also count with a subset $Z \subset Y$ where the problem
$$
    \text{minimize } f(x), \text{ subject to } x \in Z
$$
has been completely solved---in the sense that we know that the minimum $\mu_Z$ exists and we have completely determined the set of $M_Z \subset Z$ where it is attained---and we search for a technique to show that the original problem reduces to this one.

In our concrete application $X$ is the space $\mathcal{K}_0(\R^2)$ of convex bodies in the plane containing the origin as an interior point, $f$ is the area functional, $Y$ is the set of unavoidable convex polygons, and $Z \subset Y$ is the set of integer polygons that are unavoidable and convex.

Let $Y_\alpha$ be a family of subsets of $Y$ indexed by a partially ordered set $(\mathsf{A},\preceq)$ such that
\begin{enumerate}
    \item[(a)] $Y = \cup_{\alpha \in \mathsf{A}} Y_\alpha$,
    \item[(b)] for every $\alpha \in \mathsf{A}$ the set
    $\{\alpha' \in \mathsf{A}: \alpha' \prec \alpha \text{ and } Y_\alpha \not \subset Z \}$ is finite.
\end{enumerate}

\begin{proposition}\label{algorithm}
    Assume there exists a set-valued map $F : Y \rightarrow 2^Y$ such that
    \begin{enumerate}
        \item $F(y) = \emptyset$ if and only if $y \in Z$,
        \item if $y' \in F(y)$, then $f(y') \leq f(y)$,
        \item if $y \in Y_\alpha$ and $y' \in F(y)$, then $y' \in Y_{\alpha'}$ with $\alpha' \prec \alpha$.
    \end{enumerate}
    Then the minimum of $f(x)$ subject to $x \in Y$ exists and equals $\mu_Z$. If furthermore $F(y) \cap M_Z \neq \emptyset$ implies that $\mu_Z < f(y)$, then the set of minima for $f$ over $Y$ equals $M_Z$.
\end{proposition}

\begin{proof}
    Note that the map $F$ generates an algorithm in the sense that if $y \in Y$, a sequence $\{y_n\}$ satisfying
    $$
        y_0 = y \text{ and } y_n \in F(y_{n-1}) \text{ for } 1 \leq n
    $$
    is finite and ends in an element of $Z$. This follows immediately from Property~(b) of the index set $\mathsf{A}$ together with Conditions~(1) and~(3). The inequality $f(y_{n}) \leq f(y_{n-1})$ in Condition (2) implies that $\mu_Z \leq f(y)$ for all $y \in Y$ and, therefore, the restriction of $f$ to $Y$ attains $\mu_Z$ as its minimum value.

    In order to prove the last statement in the proposition, note that if the sequence starting with $y$ does not end in an element of $M_Z$, then $\mu_Z < f(y)$. If $y \not \in M_Z$ and the sequence $\{y_n\}_{n = 0}^N$ ends in $M_Z$, then $y_N \in F(y_{N-1}) \cap M_Z$ and $\mu_Z < f(y_{N-1}) \leq f(y)$.
\end{proof}

The interesting idea that optimization algorithms can be expressed in terms of the dynamics of set-valued maps is due to W.I.~Zangwill (see~\cite{Zangwill:1969, Zangwill:book}) and has proved very useful in nonlinear programming. However, we ask the reader not to see in the preceding arguments anything beyond a convenient organizational tool for our proof of Theorem~\ref{thm:I} and a justification of its description as an \term{algorithm}.

\subsubsection{Decomposition of the space of unavoidable convex polygons}
As outlined in the preceding paragraphs, our first task in the description of the algorithm is to decompose the space $Y$ of unavoidable convex polygons as a union of subsets $Y_\alpha$ indexed by a partially ordered set $(\mathsf{A},\preceq)$ in such a way that if $Z$ denotes the set of unavoidable convex polygons with integer vertices, then for every $\alpha \in \mathsf{A}$ the set
$$
    \{\alpha' \in \mathsf{A}: \alpha' \prec \alpha \text{ and } Y_\alpha \not \subset Z \}
$$
is finite. With this aim we propose the following

\begin{definition}\label{weight-type}
    Given a planar polygon $P$ and a vertex $x \in P$, we define the \term{weight} of $x$ in $P$, denoted by $\weight(x,P)$, as the number of integer lines supporting $P$ at $x$. We shall say that $P$ has \term{type} $\alpha = (n,m,k)$ if it has $n$ vertices, $m$ of which are non-integer points, and the sum of the weights of these non-integer vertices is $k$. By convention, if $m = 0$, we shall set $k = 0$. Finally we denote by $Y_\alpha=Y_{(n,m,k)}$ the set of unavoidable convex polygons of type $\alpha=(n,m,k)$.
\end{definition}

Note that the set of unavoidable convex polygons with integer vertices is decomposed as
$$
    Z = \bigcup_{n \geq 3} Y_{(n,0,0)} .
$$

If $(n', m', k')$ and $(n,m,k)$ are two types, we shall say that $(n', m', k')$ is strictly smaller than $(n,m,k)$ and write $(n', m', k') \prec (n,m,k)$ if $n' < n$, or $n' = n$ and $m' < m$, or $n' = n$, $m' = m$, and $ k' > k$. Remark that this is just the lexicographic order on the set of triplets $(n,m,-k)$.

The following lemma shows that a given type is strictly greater to only a finite number of other types that do not correspond to integer polygons.

\begin{lemma}\label{weight-vertex}
    Let $P \subset \R^2$ be an unavoidable convex polygon. If two or more distinct integer lines support $P$ at a vertex $x$, then $x$ is an integer point.
\end{lemma}

\begin{proof}
    Integer lines supporting $P$ at the vertex $x$ are dual to integer points contained in the edge of $P^*$ dual to $x$. Let us denote by $u$ and $v$ two of these points such that no other integer point lies between them. Since the polygon $P$ is unavoidable, the origin is the only integer point in the interior of $P^*$. It follows that the only integer points in the closed triangle formed by the origin, $u$, and $v$ are precisely these three vertices. This implies that $u$ and $v$ form a basis of $\Z^2$ (seen here as the dual lattice) and that $\det(u,v) = 1$. The vertex $x$ is then the unique solution of the system of equations $u \cdot x = v \cdot x = 1$ and so it must be an integer point.
\end{proof}

From this lemma we conclude that the set of types of unavoidable convex polygons is given by triples of non-negative integers $n$, $m$, and $k$ such that
$$
    3 \leq n, \ m \leq n, \ 0 \leq k \leq m.
$$

\subsubsection{Deformations of unavoidable polygons}
Having found the previous decomposition of the set $Y$ of unavoidable convex polygons, it remains for us to find the set-valued map $F : Y \rightarrow 2^Y$ satisfying the hypotheses of Proposition~\ref{algorithm}. The general idea is that if $P$ is not an integer polygon, each polygon in $F(P)$ will be obtained from $P$ by a deformation that depends on a non-integer vertex $x \in P$ and a line $\ell$ supporting $P$ at $x$.

Let us start by describing these deformations. At the moment we do not need to assume that the polygons in question are unavoidable or even that they contain the origin as an interior point.

Let $P$ be a convex, planar polygon. We denote its set of vertices by
$$
    \vertex(P) := \{x_0,\ldots,x_k\}
$$
and assume they are enumerated in positive, cyclic order: if $x_i$ is any vertex, the vectors $x_{i-1} - x_i$ and $x_{i + 1} - x_i$ form a positive basis of $\R^2$. Note that because of the cyclic order $x_{k+1} = x_0$ and hence the neighbors of $x_0$ are $x_1$ and $x_k$.

Given a vertex, which we may take to be $x_0$, we shall define a \term{supporting vector} to be a unit vector $v$ such that the line $\ell(x_0,v) := \{x_0 + tv : t \in \R \}$ supports $P$. Our deformations will be obtained by sliding the vertex $x_0$ along this line. To be precise, we define $P_t$ as the convex hull of the points $x_0 + tv$, $x_1,\ldots, x_k$:
$$
    P_t := \ch \{x_0 + tv,x_1,\ldots,x_k\} .
$$

Although we cannot restrict ourselves to small deformations, we shall restrict the range of the parameter $t$ to the interval $[\tau_{-},\tau_{+}]$, where

\begin{tabular}{ r l }
    $\tau_{-}$ & $:= \inf \{t : \vertex(P_s) = \{x_0 + sv,x_1,\ldots,x_k\} \text{ for } s \in (t,0] \}$ and \\
    $\tau_{+}$ & $:= \sup \{t : \vertex(P_s) = \{x_0 + sv,x_1,\ldots,x_k\} \text{ for } s \in [0,t) \}.$ \\
\end{tabular}

The following remarks are trivial, but they are nevertheless useful to keep in mind:
\begin{itemize}
    \item $\tau_{-} < 0$ and $\tau_{+} > 0$;
    \item either one of $\tau_{-}$ and $\tau_{+}$ can be infinite, however the only case in which both are infinite is when $P$ is a triangle and the support line $\ell(x_0,v)$ is parallel to the side opposite to the vertex $x_0$;
    \item if $\tau_+$ (resp. $\tau_-$) is finite, the polygon $P_{\tau_+}$ (resp. $P_{\tau_-}$) has exactly one vertex less than $P$ and contains all the vertices of $P$ with the exception of $x_0$.
\end{itemize}

\medskip
\begin{center}
    \includegraphics{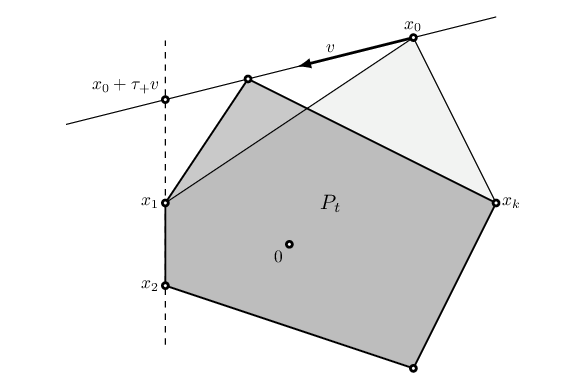}
    \captionof{figure}{Deforming a polygon.}
\end{center}
\medskip

As long as we restrict $t$ to the interval $[\tau_{-},\tau_{+}]$ we have good control on the area of $P_t$:

\begin{lemma}\label{area-function}
    The restriction of $t \mapsto |P_t|$ to the interval $[\tau_{-},\tau_{+}]$ is an affine function of $t$ which is constant if and only if the support line $\ell(x_0,v)$ is parallel to the line containing $x_1$ and $x_k$.
\end{lemma}

\begin{proof}
    Decompose $P_t$ into the union of the convex hull of $\{x_1,x_2,\ldots,x_k\}$ and the triangle formed by the vertices $x_0 + tv$, $x_1$ and $x_k$. The interiors of these sets are disjoint when $t \in [\tau_{-},\tau_{+}]$ and the result follows from the formula
    $$
        |P_t| = |\ch \{x_1,x_2,\ldots,x_k\}| + \frac{1}{2} \det(x_1 -(x_0 + tv),x_k -(x_0 + tv)) .
    $$
\end{proof}

Note that if the support line $\ell(x_0,v)$ intersects the line joining the vertices $x_1$ and $x_k$, the area of $P_t$ decreases precisely when we slide the vertex $x_0 + tv$ towards the point of intersection. From now on, we will assume that the supporting vector $v$ points in the direction of this intersection, if it exists. This will ensure that the area of $P_{t}$ will not exceed that of $P$ for $ 0 \leq t < \tau_+$.

We have reached a point where we need to be more precise about the deformations that will be performed on an unavoidable polygon $P$. For convenience, let us define a \term{virtual deformation} of $P$ to be a pair $(x,v)$, where $x$ is a non-integer vertex of $P$ and $v$ is a unit vector satisfying the following conditions:
\begin{enumerate}
    \item if $\weight(x,P) = 0$, then $v$ is a support vector directed so that the area of $P_t$ does not increase along the deformation (see the remark after Lemma~\ref{area-function});
    \item if $\weight(x,P) = 1$, then $\ell(x,v)$ will be the unique integer line that supports $P$ at $x$ and the orientation of $v$ is again such that the area of $P_t$ does not increase along the deformation.
\end{enumerate}

Notice that the set of virtual deformations is empty if and only if $P$ is an integer polytope. In the following lemma and in the rest of the section we shall assume that all deformations are of the form
$$
    P_t := \ch \{x_0 + tv,x_1,\ldots,x_k\} ,
$$
where $(x_0,v)$ is a virtual deformation of $P$ and $0 \leq t \leq \tau_+ = \tau_+(x_0,v)$.

\begin{lemma}\label{stability}
    Given a virtual deformation $(x,v)$, there exists $\epsilon > 0$ such that for every $t$ in the interval $[0 , \epsilon)$ we have that
    \begin{enumerate}
        \item the deformed polygon $P_t$ is unavoidable;
        \item $\weight(x_0 + tv, P_t) = \weight(x_0,P)$.
    \end{enumerate}
\end{lemma}

\begin{proof}
    We shall assume that weight of the vertex $x_0$ is $1$ and let $\ell_1$ denote the unique integer line that supports $P$ at $x_0$. The proof in the case that $\weight(x_0,P) = 0$ is similar.

    Consider a closed disc $B$ centered at the origin and contained in the interior of $P$. Let $\epsilon_1 > 0$ be such that the interiors of polygons $P_t$ contain $B$ for $0 \leq t < \epsilon_1$. Let $\ell_1,\ldots,\ell_N$ be the integer lines (finite in number) that intersect neither $B$ nor the relative interior of the segment joining $x_1$ to $x_k$. It is easy to see that if $\epsilon$ is any a positive number less than $\epsilon_1$ and less than the distance from $x_0$ to any of the lines $\ell_2, \ldots, \ell_N$, then $P_t$ is unavoidable and $\weight(x_0 + tv, P_t) = \weight(x_0,P)$ for all $0 \leq t < \epsilon$.
\end{proof}

\subsubsection{The set-valued map $F$ and the proof of Theorem~\ref{thm:I}}
Let $P$ be an unavoidable convex polygon, let $(x_0,v)$ be a virtual deformation of $P$, and let $P_t $ $(0 \leq t \leq \tau_+)$ be the associated deformation.

Let us start by defining $T$ to be the supremum of the set of numbers $t$ $(0 \leq t \leq \tau_+)$ such that
\begin{itemize}
    \item $P_s$ is unavoidable for all $s \in [0,t)$;
    \item $\weight(x_0 + sv, P_s) = \weight(x_0,P)$ for all $s \in [0,t)$.
\end{itemize}

\begin{claim}
    $T$ is finite and strictly positive.
\end{claim}

\begin{claimproof}
    The strict positivity of $T$ follows from Lemma~\ref{stability}. In order to see that $T$ is finite we reason as follows: if $T$ is infinite we must have that $\tau_+$ is infinite and that the line defined by the vertices $x_1$ and $x_k$ is parallel to $v$. This implies that $P$ is a triangle. However, as we move the vertex $x_0 + sv$ to infinity, the origin must leave the triangle. By Lemma~\ref{origin-inside} this contradicts the fact that $P_s$ is unavoidable for all values of $s < T$.
\end{claimproof}

We now define $F(P)$ as the set of all the polytopes $P_T$ that are obtained from the previous construction starting from a virtual deformation of $P$. Since the condition of being unavoidable is closed, it follows that each $P_T$ is an unavoidable polygon and $F$ is indeed a map from $Y$ to $2^Y$.

Note that $F$ assigns the empty set to an unavoidable convex polygon if it has no virtual deformations and this happens if and only of it is an integer polygon. Moreover, by construction, the area of $P_T \in F(P)$ does not exceed that of $P$. In other words, the first two conditions of Proposition~\ref{algorithm} are fulfilled.

We must now show that the type (see~Definition~\ref{weight-type}) of $P_T$ is strictly smaller than the type of $P$. Since, by construction, we have that $P_T$ cannot have more vertices than $P$, this is an immediate consequence of the following claim.

\begin{claim}
    Either $\weight(x_0 + Tv, P_T) > \weight(x_0,P)$ or the number of vertices of $P_T$ is strictly less than that of $P$. Moreover, $\weight(x,P_T) \geq \weight(x,P)$ for every vertex $x$ common to $P$ and $P_T$.
\end{claim}

\begin{claimproof}
    If $T = \tau_+$, then $P_T$ has one vertex less than $P$. If $T$ is strictly less than $\tau_+$ and $\weight(x_0 + Tv, P_T) = \weight(x_0,P)$ then, Lemma~\ref{stability} tells us that there exists an $\epsilon > 0$ such that $P_{T+ s}$ is unavoidable and
    $$\
        \weight(x_0 + (T + s)v, P_{T+ s}) = \weight(x_0 P) \text{\quad for all } 0 \leq s < \epsilon .
    $$
    This contradicts the definition of $T$ and thus we conclude that
    $$
        \weight(x_0 + Tv; P_T) > \weight(x_0,P).
    $$

    To see that $\weight(x,P_T) \geq \weight(x,P)$ for every vertex $x$ common to $P$ and $P_T$ first note that $\weight(x,P_t) \geq \weight(x,P)$ for $0 \leq t < T$. Indeed, the only way for the weight of $x$ to diminish during the deformation is that a fixed integer line $\ell$ that supports the polygon at $x$ at some instant $t$ intersects its interior at some later time $t' < T$. By continuity, this cannot happen without increasing the weight of the only vertex that moves during the deformation. However, this increase is explicitly forbidden by the definition of $T$. In order to conclude that $\weight(x,P_T) \geq \weight(x,P)$ we just need to remark that a line that supports $P_t$ at $x$ for all $t \in [0, T)$ must also support $P_T$ at $x$.
\end{claimproof}

Proposition~\ref{algorithm} now tells us that the minimal area of an arbitrary unavoidable convex polygon equals the minimal area of an unavoidable, convex, integer polygon. By Lemma~\ref{integer-polytope}, this minimal area is $3/2$ and, therefore, \assertion{the area of an unavoidable convex polygon is at least $3/2$. }

We conclude the proof of Theorem~\ref{thm:I} by characterizing the set of unavoidable convex polygons of minimal area. By the second part of Proposition~\ref{algorithm}, we just need to verify that if $F(P)$ contains a basic triangle, then the area of $P$ is strictly greater than $3/2$.

By applying a unimodular transformation if necessary, we may assume that the basic triangle $P_T$ has vertices $(1,0)$, $(0,1)$, and $(-1,-1)$. Moreover, since the subgroup of $GL(2,\Z)$ that fixes this triangle acts transitively on its vertices, we may assume that vertices $(1,0)$ and $(0,1)$ are vertices of $P$ that remained fixed along the deformation while the vertex $(-1,-1)$ is the vertex $x_0 + Tv$.

The key remark is that the points $(1,0)$ and $(0,1)$ and any other vertices of $P$ different from $x_0$ must also belong to the triangle $P_T$. Because of this, if $P$ contains $(-1,-1)$, then its area is strictly greater than $3/2$. Likewise, if $(-1,-1)$ is not in $P$, then the vertex $x_0$ must lie in the region $x < -1$, $y < -1$ or it would not be possible for the polygon $P$ to intersect the lines $x = -1$ and $y = - 1$, which it must since it is unavoidable. In this case, we also have that the area of $P$ is strictly greater than $3/2$.
\qed

\section{Systolic geometry of optical hypersurfaces}\label{systolic}

In this section we show that Theorem~\ref{thm:I} extends to a sharp isosystolic inequality for optical hypersurfaces in the cotangent bundle of the two-dimensional torus and that much of systolic geometry can be likewise extended to the setting of optical hypersurfaces in cotangent bundles of compact manifolds. In the equivalent language of Finsler metrics, we shall prove the Finsler analogue of Loewner's isosystolic inequality and show that much of systolic geometry can be extended to the setting of Finsler metrics \term{that are not necessarily reversible}.

\subsection{Optical hypersurfaces and Finsler metrics} \label{duality-Finsler-optical}

Let us start by recalling some basic definitions and constructions.

\begin{definition}
    A smooth hypersurface $\Sigma \subset T^*M$ in the cotangent bundle of a smooth $n$-dimensional manifold $M$ is said to be an \term{optical hypersurface} if each intersection $\Sigma_x := \Sigma \cap T_x^* M$ ($x \in M$) is a quadratically convex hypersurface in $T_x^*M$ enclosing the origin. If $\Sigma_x \subset T_x^* M$ is $0$-symmetric for every $x \in M$, we shall say that the hypersurface $\Sigma$ is \term{reversible}.
\end{definition}

\begin{definition}\label{Finsler-metric}
    A \term{Finsler metric} on a manifold $M$ is a continuous function $F : TM \rightarrow [0,\infty)$ that is smooth outside the zero section and whose restriction to every tangent space is a Minkowski norm. If $F(-v_x) = F(v_x)$ for all tangent vectors, we shall say that $F$ is a \term{reversible Finsler metric}.
\end{definition}

\noindent
\textit{Remark.} It is often useful to also consider \term{continuous Finsler metrics}: continuous functions on $TM$ whose restriction to every tangent space is a norm.

The notion of duality allows us to pass back and forth from optical hypersurfaces to Finsler metrics and from reversible optical hypersurfaces to reversible Finsler metrics:

Given an optical hypersurface $\Sigma \subset T^*M$ and a point $x \in M$, we may consider $\Sigma_x \subset T_x^* M$ as the unit sphere of a Minkowski norm $\|\cdot \|^*_x$ on $T_x^*M$. If we denote the dual norm on $T_x M$ by $\| \cdot \|_x$ and define $F : TM \rightarrow [0,\infty)$ by $F(v_x) := \|v_x\|_x $, we obtain a Finsler metric. This Finsler metric and the Hamiltonian $H_{{}_\Sigma}(p_x) = \|p_x\|^*_x$ are said to be dual.

If we start with a Finsler metric $F$, the first step in constructing an optical hypersurface is to consider the unit disc bundle
$$
    \D(F) := \{v \in TM: F(v) \leq 1 \} ,
$$
which is the disjoint union of the convex bodies
$$
    \D_x(F) := \D(F) \cap T_xM \ \ (x \in M).
$$
The \term{unit co-disc bundle} $\D^*(F)$ is defined as the union of the dual bodies $\D_x^*(F) \subset T^*_x M$ and its boundary, the \term{unit co-sphere bundle} $\Sph^*(F) := \partial \D^*(F)$, is an optical hypersurface in the cotangent bundle of $M$.

The description of a Riemannian or Finsler metric in terms of its unit co-sphere bundle provides us with useful characterizations of its geodesic spray and its volume.

The pullback of the canonical $1$-form $\alpha$ on $T^* M$ to the optical hypersurface $\Sigma$ is a contact form, which we denote again by $\alpha$. If $\Sigma$ is the unit co-sphere bundle of the Finsler metric $F$, the orbits of the Reeb vector field $R_\alpha : \Sigma \rightarrow T\Sigma$---defined by the equations $d\alpha(R_\alpha,\cdot) = 0$ and $\alpha(R_\alpha) = 1$---project down to geodesics in $(M,F)$ that are parameterized by arc length. In particular, the action of a periodic Reeb orbit $\gamma$---defined as the integral of the contact form $\alpha$ over $\gamma$---equals the length or period of the corresponding closed geodesic.

The volume of an optical hypersurface $\Sigma$ of dimension $2n-1$ is defined as the integral of the volume form $\alpha \wedge (d\alpha)^{n-1}$. By Stokes's formula, this is the volume of the open subset of $T^*M$ enclosed by $\Sigma$ measured with the symplectic volume form $\omega^n$. We warn the reader that this is $n!$ times what is usually taken as the symplectic volume form.

If $\Sigma$ is the unit co-sphere bundle of the Finsler metric $F$, its volume is $n! \varepsilon_n \vol(M,F)$, where $\varepsilon_n$ is the volume of the Euclidean unit sphere of dimension $n$ and $\vol(M,F)$ is the \term{Holmes-Thompson volume} of $(M,F)$. For Riemannian metrics this notion of volume coincides with the Hausdorff measure, but this is no longer so for more general reversible Finsler metrics. In that case we have the following useful result of Dur\'an~\cite{Duran:1998}:

\begin{theorem}[Dur\'an]
    The Holmes-Thompson volume of a reversible Finsler manifold does not exceed its Hausdorff measure. When the Hausdorff measure is finite, it equals the Holmes-Thompson volume if and only if the manifold is Riemannian.
\end{theorem}

Remark that there is no established definition for the Hausdorff measure of an asymmetric metric space and the comparison between Holmes-Thompson volume and Hausdorff measure can only be made for reversible Finsler metrics.

\subsection{The geometry of numbers as the geometry of Finsler tori} \label{Loewner-subsection}

A simple and far-reaching observation that the arithmetic study of positive-definite quadratic forms can be seen as the study of flat Riemannian tori. The full meaning of this metaphor is more striking if, like Minkowski, we think of ellipsoids as particular examples of convex bodies. Consider, for instance, the following Finsler-geometric translation (see Gromov~\cite[p.~295]{Gromov:1996}) of Minkowski's lattice-point theorem: \assertion{if the Hausdorff measure of a flat and reversible $n$-dimensional Finsler torus equals $\varepsilon_n/2^n$, it carries a non-contractible periodic geodesic of length at most $1$}.

\begin{proof}[Translation]
    Let us write the flat, reversible Finsler torus as the quotient of a reversible normed space $(\R^n,\|\cdot\|)$ by a lattice $\Lambda$ and remark that:
    \begin{itemize}
        \item The Hausdorff measure of $(\R^n,\|\cdot\|)/\Lambda$ equals $\varepsilon_n\det(\Lambda)/|B|$, where $B$ is the unit ball in $(\R^n,\|\cdot\|)$.
        \item All periodic geodesics in $(\R^n,\|\cdot\|)/\Lambda$ are non-contractible and correspond to non-zero lattice points. The lengths of these geodesics equal the norms of the corresponding points.
    \end{itemize}

    Since, by hypothesis, $|B|/\det(\Lambda)= 2^{n}$, it follows from Minkowski's lattice-point theorem that $B$ contains a non-zero lattice point. This point corresponds to a non-contractible periodic geodesic of length at most $1$.
\end{proof}

Equally simple translations of Theorems~I and~II yield the following results:

\begin{itemize}
    \item \assertion{If the Holmes-Thompson area of a flat two-dimensional Finsler torus equals $3/2\pi$, it carries a non-contractible periodic geodesic of length at most $1$.}
    \item \assertion{If the Holmes-Thompson volume of a flat $n$-dimensional Finsler torus equals $(\pi/4)^n /\varepsilon_n n!$, it carries a non-contractible periodic geodesic of length at most $1$.}
\end{itemize}

It is natural to ask whether these results extend to Finsler tori that are not flat. Remarkably, the work of D.~Burago and S.~Ivanov on the minimality of flats in normed spaces (see~\cite{Burago-Ivanov:2002}) enables us to show that this is the case in two dimensions.

\begin{maintheorem}{IV}\label{thm:IV}
    If an optical hypersurface in the cotangent bundle of the $2$-torus bounds a volume $V$, it carries a periodic characteristic whose action is at most $\sqrt{V/3}$. Equivalently, if the Holmes-Thompson area of a Finsler $2$-torus equals $3/2\pi$, it carries a periodic geodesic of length at most $1$.
\end{maintheorem}

The key idea in the proof is to reduce the problem to the flat case by a standard homogenization technique which we now review for the reader's convenience. Details of the constructions that follow can be found in Gromov~\cite[pp.~245--261]{Gromov:1999}.

Given an optical hypersurface $\Sigma$ in the cotangent bundle of a compact manifold $M$, consider the $1$-homogeneous function $H_{{}_\Sigma} : T^*M \rightarrow [0,\infty)$ that is constantly equal to $1$ on $\Sigma$ and define a norm $\|\cdot\|_{\rm st}^*$ on $\mathsf{H}^1(M;\R) $, the first real cohomology of $M$, by
$$
    \|{\tt a}\|_{\rm st}^* := \inf \{ \max_{x \in M} H_{{}_\Sigma}(\eta(x)) : \eta \text{ is a closed $1$-form with } [\eta] = {\tt a} \} .
$$
The dual norm $\|\cdot\|_{\rm st} : \mathsf{H}_1 (M;\R) \rightarrow [0,\infty)$ is the \term{stable norm} of the Finsler metric dual to $H_{{}_\Sigma}$.

On the $n$-dimensional torus $\T^n$ we can identify real cohomology classes of degree $1$ with translation invariant $1$-forms and the preceding construction can be interpreted as a homogenization of the Hamiltonian $H_{{}_\Sigma}$. A more direct description of the homogenized Hamiltonian can be given if we note that every closed $1$-form on the torus is the sum of an invariant $1$-form and the differential of a function:

\begin{definition}
    If $\Sigma$ is an optical hypersurface in the cotangent space of the torus $\T^n$ and $H_{{}_\Sigma} : T^*\T^n \rightarrow [0,\infty)$ is the $1$-homogeneous function that is constantly equal to $1$ on $\Sigma$, we define the \term{homogenized Hamiltonian}
    $$
        \widehat{H}_{{}_\Sigma}(x,p) := \inf_{f \in C^1(\T^n)} \max_{x \in \T^n} H_{{}_\Sigma}(x,p + df(x)) .
    $$
\end{definition}

The Hamiltonian $\widehat{H}_{{}_\Sigma}$ is independent of the variable $x$ and, as was mentioned before, can be considered as a norm in $\R^{n*} = \mathsf{H}^1(M;\R)$. However, we prefer to consider $\widehat{H}_{{}_\Sigma}$ as the Hamiltonian of a flat (non-smooth) Finsler metric on $\T^n$.

In this terminology, the main result of Burago and Ivanov~\cite{Burago-Ivanov:2002} can be stated as follows:

\begin{theorem}[Burago and Ivanov]\label{Burago-Ivanov-HT}
    The volume enclosed by an optical hypersurface $\Sigma$ in $T^*\T^2$ is no less than the volume enclosed by the energy surface $\widehat{H} _{{}_\Sigma} = 1$ for the homogenized Hamiltonian $\widehat{H} _{{}_\Sigma}$. Equivalently, the Holmes-Thompson area of a Finsler metric on the $2$-dimensional torus is no less than the Holmes-Thompson area of the flat metric defined by its stable norm.
\end{theorem}

To be precise, Burago and Ivanov state this theorem for reversible Finsler metrics, but a detailed study of their proof reveals that it extends unchanged to the non-reversible case.

\begin{proof}[Proof of Theorem~\ref{thm:IV}]
    A classic result in systolic geometry (see Lemma~4.32 in~\cite[p.~260]{Gromov:1999}) states that the length of the shortest non-contractible periodic geodesic in a $2$-dimensional Finsler torus equals the length of the shortest non-contractible periodic geodesic for the flat metric defined by its stable norm (i.e., the stable systole). Together with Theorem~\ref{Burago-Ivanov-HT}, this immediately implies that if $(\T^2, F)$ is a Finsler torus we can always find a flat torus $(\T^2,\widehat{F})$ whose area does not exceed that of $(\T^2,F)$ and such that the shortest non-contractible geodesics in both tori have the same length. In other words, we have reduced the proof of the theorem to the case of flat Finsler tori and this is precisely the content of Theorem~\ref{thm:I}.

    The equivalent statement for optical hypersurfaces in the cotangent bundle of the torus follows from their relation to Finsler metrics explained in~\ref{duality-Finsler-optical}.
\end{proof}

The idea of the preceding proof is taken from Sabourau~\cite{Sabourau:2010}, where it is used to prove the following sharp systolic inequality for reversible Finsler metrics on $\T^2$:

\begin{theorem}[Sabourau]\label{Sabourau}
    If a reversible optical hypersurface in the cotangent bundle of the $2$-torus bounds a volume $V$, it carries a periodic characteristic whose action is at most $\sqrt{V}/2$. Equivalently, if the Holmes-Thompson area of a Finsler $2$-torus equals $2/\pi$, it carries a periodic geodesic of length at most $1$.
\end{theorem}

The inequalities in Theorem~\ref{thm:IV} and Theorem~\ref{Sabourau} are sharp. In the first case, in agreement with Theorem~\ref{thm:I}, equality is attained for the flat Finsler metric in $\R^2/\Z^2$ induced by the norm whose dual unit ball in $\R^{2*}$ is the triangle with vertices $(1,0)$, $(0,1)$, and $(-1,-1)$. In the second case, equality is attained for the flat Finsler metric induced by the norm whose dual unit ball is the square with vertices $\pm(1,0)$, $\pm(0,1)$. Note that this metric can also be described as the flat Finsler metric on $\R^2/\Z^2$ induced by the $\ell_\infty$ norm on the plane, or the Finsler metric dual to the Hamiltonian
$$
    H(x_1,x_2,p_1,p_2) = |p_1| + |p_2| .
$$
This Hamiltonian is smooth on a dense open subset of $T^*\T^2$ and on this set Hamilton's equations of motion define a periodic flow in which all orbits have the same minimal period. For the reader who is acquainted with the work of \'Alvarez Paiva and Balacheff (\cite{Alvarez-Balacheff:2011}) on the special role of Zoll Finsler manifolds and regular contact manifolds in systolic geometry, we shall now explain why this remark plays an important heuristic role in our investigations.

\subsection{* Zoll tori}\label{Zoll-tori-subsection}

In~\cite{Alvarez-Balacheff:2011} \'Alvarez Paiva and Balacheff prove that if a (smooth) Finsler metric $F$ on a compact manifold $M$ has minimal volume among all the Finsler metrics on $M$ for which the length of the shortest periodic geodesic equals some fixed constant, then $(M,F)$ is a \term{Zoll Finsler manifold}: all its geodesics are periodic with the same minimal period. Since most manifolds do not admit smooth Zoll metrics, this means that on most manifolds extremal Finsler metrics, if they exist, are not smooth. However, as illustrated by the remarks at the end of~\S\ref{Loewner-subsection}, extremal metrics may still be Zoll metrics if we slightly broaden the definition:

\begin{definition}\label{generalized-Zoll}
    A continuous Finsler metric on a compact manifold $M$ is a \term{generalized Zoll metric} if its dual Hamiltonian $H: T^*M \rightarrow [0,\infty)$ is smooth on an open set of full measure $\mathcal{O} \subset T^*M$ on which Hamilton's equations of motion define a periodic flow and all orbits have the same minimal period.
\end{definition}

We do not know whether every compact manifold admits a generalized Zoll Finsler metric, but it is quite easy to construct such metrics on tori:

\begin{proposition}\label{flat-Zoll-tori}
    Consider a norm $\|\cdot\| : \R^n \rightarrow [0,\infty)$ with unit ball $B$. The induced flat Finsler metric on the torus $\T^n = \R^n/\Z^n$ is a generalized Zoll metric of period $1$ if and only if $B$ is a polytope all of whose vertices are primitive elements of the integer lattice $\Z^n$.
\end{proposition}

\begin{proof}
    Note that the Hamiltonian $H : \T^n \times \R^{n*} \rightarrow [0,\infty) $ dual to the flat Finsler metric is simply $H(x,p) = \|p\|^*$. Since the dual norm is homogeneous and convex, its differential is defined on a conic set of full measure $U \subset \R^{n*}$. On $\T^n \times U$ Hamilton's equations,
    $$
        \dot{x} = \frac{\partial H(x,p)}{\partial p} \text{\quad and\quad} \dot{p} = - \frac{\partial H(x,p)}{\partial x} = 0 ,
    $$
    define the flow $\varphi_t(x,p) = (x + t \partial H(x,p)/\partial p, p)$. This flow is periodic and all orbits have the same minimal period equal to $1$ if and only if $\partial H(x,p)/\partial p$ is a primitive element of $\Z^n$ for every initial condition $(x,p) \in \T^n \times U$.

    The homogeneity of $H(x,p) = \| p \|^*$ allows us to restrict our attention to those initial conditions $(x,p)$ with $H(x,p) = \|p\|^* = 1$ (i.e., to the hypersurface $\T^n \times \partial B^*$). In this case, the differential $\partial H(x,p)/\partial p$ is the point in $\partial B$ corresponding to the unique hyperplane supporting $B^*$ at $p$.

    If $B$ is a polytope, by duality, its vertices correspond to the facets of the dual polytope $B^*$. Moreover, if $p$ is in the relative interior of a facet of $B^*$, then $\|\cdot\|^*$ is smooth at $p$ and its differential is precisely the vertex of $B$ corresponding to the facet. Thus, when $B$ is a polytope, the Hamiltonian $H$ is smooth on an open set of full measure, and the flat Finsler metric defined by the norm $\|\cdot\|$ is a generalized Zoll metric of period $1$ if and only if the vertices of $B$ are primitive elements of $\Z^n$.

    It remains for us to show that if at every point of $U$ the differential of the norm $\|\cdot\|^*$ is a primitive element of $\Z^n$, then $B$ must be a polytope. As before, we may restrict our attention to $U \cap \partial B^*$: those points on the boundary of $B^*$ at which there is a unique support plane. If $p$ is any of these points, $\partial H(x,p)/\partial p$ lies on the (compact) boundary of $B$. If $\partial H(x,p)/\partial p$ also takes values on the discrete set of primitive element of $\Z^n$, then it can only take a finite set of values. This implies that a finite number of hyperplanes support $B^*$ at a dense set of points. Since the set of points at which a given hyperplane supports a convex body is closed, we conclude that a finite set of hyperplanes support $B^*$ at every point of its boundary and, therefore, $B^*$ and $B$ must be polytopes.
\end{proof}

At an informal level, this proposition and the identification in~\cite{Alvarez-Balacheff:2011} of Zoll manifolds as the only smooth extremal Finsler metrics form the basis of the present work: in Lemma~\ref{integer-polytope} we established the sharp systolic inequality for flat generalized Zoll metrics on the torus $\R^n/\Z^n$, in the algorithmic proof of Theorem~\ref{thm:I} we showed that any flat metric on the $2$-torus can be deformed to a flat generalized Zoll metric without decreasing the systolic volume, and, in Theorem~\ref{thm:IV} we showed that the systolic volume of any Finsler metric in the $2$-torus is no less than the systolic volume of an associated flat metric. In other words, we verified that the $2$-torus admits extremal metrics that are Zoll---in the generalized sense of Definition~\ref{generalized-Zoll}---and used this to determine the sharp isosystolic inequality for Finsler $2$-tori.

\subsection{Systolic capacities and coarse isosystolic inequalities}\label{coarse-inequalities-subsection}

The functional that assigns to every compact Riemannian manifold $(M,g)$ the length of its shortest periodic geodesic is merely lower semi-continuous as a function of the metric and is not amenable to a very detailed analysis (nevertheless, see~\cite{Alvarez-Balacheff:2011}). Its main shortcoming from the point of view of systolic geometry is its lack of monotonicity: if $g \leq g'$ are two metrics on $M$, the shortest periodic geodesic of $(M,g)$ is not necessarily shorter than the shortest periodic geodesic of $(M,g')$. For this reason, it is usual to work with Riemannian invariants that are monotone and which are closely related to the length of a periodic geodesic. We shall call such invariants \term{systolic capacities}.

\begin{definition}\label{systolic-capacity}
    A real-valued function $\kappa$ defined on the class of smooth Finsler metrics on a manifold $M$ is said to be a systolic capacity if it satisfies the following conditions:
    \begin{enumerate}
        \item homogeneity: $\kappa(M,\lambda F) = \lambda \, \kappa(M,F)$ for positive numbers $\lambda$;
        \item monotonicity: $\kappa(M,F) \leq \kappa(M,G)$ whenever $F \leq G$;
        \item representability: $\kappa(M,F)$ equals the length of a finite collection of periodic geodesics (counted with multiplicity).
    \end{enumerate}
\end{definition}

The three main examples of systolic capacities are the \term{homotopic systole}---the length of the shortest non-contractible periodic geodesic---defined on multiply-connected Finsler manifolds, the \term{diastole} defined for any compact Finsler manifold (see pages~334--335 of~\cite{Bott:1982} for a clear and elegant account of this invariant and its fundamental role in the proof of the Luysternik-Fet theorem), and the \term{diastole over $1$-cycles} $\dias_{\mathcal{Z}}(M^2,F)$ defined for compact orientable Finsler surfaces. For the delicate construction of this last invariant, which is due to F.~Almgren~\cite{Almgren:1962} and J. Pitts~\cite{Pitts:1981}, we refer the reader to Part II of~\cite{Calabi-Cao:1992} and~\cite[\S 2]{Balacheff-Sabourau:2010}.

The monotonicity and homogeneity of a systolic capacity allow us to extend its domain of definition to continuous Riemannian or Finsler metrics. Indeed, this is an immediate consequence of the following folklore result:

\begin{lemma}\label{approximation}
    Let $F$ be a continuous Finsler metric on a (smooth) manifold $M$. For every $\epsilon > 0$, there exists a smooth Finsler metric $F_\epsilon$ on $M$ such that $F_\epsilon < F < (1 + \epsilon) F_\epsilon$. Moreover, if $F$ is a continuous Riemannian metric, $F_\epsilon$ may be chosen to be a smooth Riemannian metric.
\end{lemma}

Another important consequence of the monotonicity of systolic capacities is the following

\begin{theorem}\label{inequality-reversible}
    Assume that we are given a compact $n$-dimensional manifold $M$ and a systolic capacity $\kappa$ for which there exists a constant $C > 0$ such that $ \kappa(M,\sqrt{g}) \leq C \sqrt[n]{\vol(M,\sqrt{g})}$ for every Riemannian metric $g$ on $M$, then
    $$
        \kappa(M,F) \leq \frac{C \sqrt[n]{\varepsilon_n n!} }{2} \sqrt[n]{\vol(M,F)}
    $$
    for every reversible Finsler metric $F$ on $M$.
\end{theorem}

The easy translation to the language of optical hypersurfaces yields the following

\begin{corollary}\label{corollary-reversible}
    Under the hypotheses of Theorem~\ref{inequality-reversible}, any reversible optical hypersurface $\Sigma \subset T^*M$ bounding a volume $V$ carries a periodic characteristic whose action is at most $C \sqrt[n]{V}/2$.
\end{corollary}

\begin{proof}[Proof of Theorem~\ref{inequality-reversible}]
    Recall that the \term{inner Loewner ellipsoid} of a convex body $K$ in a $n$-dimensional vector space $V$ is the unique ellipsoid $E$ contained in $K$ and minimizing the ratio $|K|/|E|$. Let us consider the map $x \mapsto g_x$ that assigns to every point $x \in M$ the quadratic form for which the ellipsoid $g_x(v,v) \leq 1$ is the inner Loewner ellipsoid of the unit ball $\D_x(F) \subset T_xM$. This map defines a continuous Riemannian metric $g$ on $M$ (for a simple proof of this statement, see S.~Ivanov's answer in \url{http://mathoverflow.net/questions/127114}). By construction, this metric is larger than $F$ (i.e., $F \leq \sqrt{g}$), and hence
    $$
        \kappa(M,F) \leq \kappa(M,\sqrt{g}) \leq C \sqrt[n]{\vol(M,\sqrt{g})} .
    $$
    We shall now compare the volume of the Riemannian manifold $(M,g)$ and the Holmes-Thompson volume of the Finsler manifold $(M,F)$ by comparing the Holmes-Thompson volume densities
    $$
        \mu_{\sqrt{g}}, \, \mu_F^{ht} : \bigwedge^n(TM) \longrightarrow [0,\infty)
    $$
    at each point $x \in M$. Following~\cite[p.~11]{Alvarez-Thompson:2004}, the density $\mu_F^{ht}$ is given by the formula
    $$
        \mu_F^{ht}(v_1\wedge v_2 \wedge \cdots \wedge v_n) =
            \varepsilon_n^{-1} \frac{|\D_x^*(F)|}{|v_1^* \wedge v_2^* \wedge \cdots \wedge v_n^*|} .
    $$
    The Holmes-Thompson volume density of a Riemannian manifold is its standard volume density and, therefore, $\mu_{\sqrt{g}}$ is given by a similar formula. In fact, recalling that $g_x(v,v) \leq 1$ is the inner Loewner ellipsoid $E_x$ of the unit ball $\D_x(F) \subset T_xM$,
    $$
        \mu_{\sqrt{g}}(v_1\wedge v_2 \wedge \cdots \wedge v_n) =
            \varepsilon_n^{-1} \frac{|E_x^*|}{|v_1^* \wedge v_2^* \wedge \cdots \wedge v_n^*|} .
    $$
    We now remark that $E_x^*$ is the \term{outer Loewner ellipsoid} of the convex body $\D_x^*(F)$ and hence, at each point $x \in M$, the function $\mu_{\sqrt{g}}/\mu_F^{ht}$ is the \term{outer volume ratio} $|E_x^*|/|\D_x^*(F)|$, which is at most $\varepsilon_n n! /2^n$ by a theorem of K.~Ball and F.~Barthe (see Theorem~5 in~\cite{Ball:2001}). Integrating the volume densities over $M$, we obtain
    $$
        \vol(M,\sqrt{g}) \leq \frac{\varepsilon_n n!}{2^n} \vol(M,F) ,
    $$
    which completes the proof of the theorem.
\end{proof}

In the preceding proof, we could have used John's ellipsoid theorem instead of the upper bound for the outer volume ratio to obtain a result valid for all definitions of volume on a reversible Finsler manifold, albeit with a worse estimate. The real advantage of using the Holmes-Thompson volume is revealed by the fact that isosystolic inequalities for reversible Finsler metrics yield isosystolic inequalities for general Finsler metrics.

\begin{theorem}\label{inequality-non-reversible}
    Assume that we are given a compact $n$-dimensional manifold $M$ and a systolic capacity $\kappa$ for which there exists a constant $C > 0$ such that $\kappa(M,F) \leq C \sqrt[n]{\vol(M,F)}$ for every reversible Finsler metric $F$ on $M$, then
    $$
        \kappa(M,G) \leq C \sqrt[n]{(2n)!/(n!)^2} \sqrt[n]{\vol(M,G)}
    $$
    for every Finsler metric $G$ on $M$.
\end{theorem}

In terms of optical hypersurfaces, we have the following

\begin{corollary}
    Under the hypotheses of Theorem~\ref{inequality-non-reversible}, any optical hypersurface $\Sigma \subset T^*M$ bounding a volume $V$ carries a periodic characteristic whose action is at most $C \sqrt[n]{(2n)!/(n!)^2} \sqrt[n]{V}$.
\end{corollary}

\begin{proof}[Proof of Theorem~\ref{inequality-non-reversible}]
    Consider the symmetrized Finsler metric defined by
    $$
        F(v) := G(v) + G(-v).
    $$
    At every $x \in M$ the dual unit ball $\D^*_x(F)$ for the reversible metric $F$ is the difference body of the dual unit ball $\D_x^*(G)$ for the metric $G$ (i.e., $\D_x^*(F) = \D_x^*(G) - \D_x^*(G)$). Applying the Rogers-Shepard inequality at every cotangent space we have that
    $$
        \vol(M,F) \leq \frac{(2n)!}{(n!)^2} \vol(M,G) .
    $$
    Using that $G \leq F$, we conclude that
    $$
        \kappa(M,G) \leq \kappa(M,F) \leq C \sqrt[n]{\vol(M,F)} \leq C \sqrt[n]{(2n)!/(n!)^2} \sqrt[n]{\vol(M,G)}.
    $$
\end{proof}

We can now present Theorems~V and~VI as simple corollaries of Riemannian results.

\begin{maintheorem}{V}\label{thm:V}
    Given an essential $n$-dimensional manifold $M$, there exists a constant $C > 0$ such that every optical hypersurface $\Sigma \subset T^* M$ bounding a volume $V$ carries a closed characteristic whose action is less than $C \sqrt[n]{V}$.
\end{maintheorem}

\begin{proof}
    Gromov's celebrated isosystolic inequality (see \cite{Gromov:1983}) states that \assertion{given an essential $n$-dimensional manifold $M$, there exists a constant $C > 0$ such that the homotopic systole of any Riemannian or reversible Finsler metric $F$ on $M$ is less than $C \sqrt[n]{\vol(M,F)}$}. The extension to non-reversible Finsler metrics and optical hypersurfaces follows at once from Theorem~\ref{inequality-non-reversible}.
\end{proof}

\begin{maintheorem}{VI}\label{thm:VI}
    There exists constant $C > 0$ such that every optical hypersurface $\Sigma \subset T^* S^2$ bounding a volume $V$ carries a closed characteristic whose action is less than $C \sqrt{V}$.
\end{maintheorem}

\begin{proof}
    The main result of Balacheff and Sabourau in~\cite{Balacheff-Sabourau:2010} states that \assertion{there exists a positive constant $C \leq 10^8$ such that every Riemannian metric $g$ on $S^2$ satisfies }
    $$
        \dias_{\mathcal{Z}}(S^2, \sqrt{g}) \leq C \sqrt{\area(S^2,\sqrt{g})} .
    $$
    Theorems~\ref{inequality-reversible} and~\ref{inequality-non-reversible} imply that this inequality remains valid for Finsler metrics if we replace $C$ by $\sqrt{3\pi} \, C$. The statement for optical hypersurfaces follows from the (by now) usual translation.
\end{proof}

It is unfortunate that in the preceding proof we cannot use the much better length-area bounds for the shortest periodic geodesic on a Riemannian $2$-sphere established by Croke~\cite{Croke:1988} and refined by R. Rotman in~\cite{Rotman:2006}. However, the monotonicity of the diastole over $1$-cycles is crucial to the proof and, in general, we do not know whether a length-volume bound for the shortest periodic geodesic established for Riemannian metrics will translate to a similar bound for Finsler metrics.

\subsection{Isosystolic inequalities and the Mahler conjecture}\label{Mahler-conjecture-subsection}

The oldest open problem in systolic geometry is to determine whether a Riemannian metric on $\RP^n$ with the same volume as the canonical metric must carry a non-contractible periodic geodesic whose length is at most $\pi$. This was proved in two dimensions by P.M.~Pu~(\cite{Pu:1952}). The existence of an upper bound for the length of the shortest non-contractible geodesic follows from Gromov's isosystolic inequality~(\cite{Gromov:1983}), but it is only recently that a good upper bound has been found. Indeed, by building on work by L.~Guth, K.~Nakamura shows in~\cite{Nakamura:2013} that a Riemannian metric on $\RP^n$ with the same volume as the canonical metric must carry a non-contractible periodic geodesic whose length is at most $\sqrt[n]{n!} \, \pi$.

Also in support of the conjectured inequality, \'Alvarez Paiva and Balacheff show in~\cite{Alvarez-Balacheff:2011} that if $g_t$ is a smooth, constant-volume deformation of the canonical metric that is not formally trivial, the length of the shortest periodic geodesic of $(\RP^n,g_t)$ attains $\pi$ as a strict local maximum at $t = 0$. We recall that a smooth, one-parameter deformation $g_t$ of Riemannian metrics is formally trivial if for every positive integer $k$ there exists an isotopy $\phi_t$ for which $\phi_t^* g_0$ and $g_t$ agree to order $k$ at $t = 0$.

The results in~\cite{Alvarez-Balacheff:2011} also suggest that it is reasonable to conjecture that any reversible Finsler metric on $\RP^n$ with the same Holmes-Thompson volume as the canonical metric must carry a non-contractible periodic geodesic whose length is at most $\pi$. Moreover the reversible optical hypersurface in the cotangent bundle of the real projective $n$-space corresponding to the canonical metric encloses a volume equal to
$$
    n! \varepsilon_n \vol(\RP^n,\textit{canonical metric}\,)=(n+1)! \varepsilon_n \varepsilon_{n+1}/2=(2\pi)^n.
$$
So in terms of optical hypersurfaces, this conjecture translates as follows

\begin{conjecture}{II}\label{conj:II}
    If a reversible optical hypersurface in the cotangent bundle of real projective $n$-space encloses a volume $V$, it carries a periodic characteristic whose action is at most $\sqrt[n]{V}/2$. Moreover, this short characteristic can be chosen so that its projection onto the base manifold is a non-contractible loop.
\end{conjecture}

In two dimensions, this conjecture has been proved by S.~Ivanov (see~\cite{Ivanov:2002} and~\cite{Ivanov:2011}). Moreover, Nakamura's result and Corollary~\ref{corollary-reversible} immediately yield the following

\begin{theorem}\label{Nakamura-ABT}
    If a reversible optical hypersurface in the cotangent bundle of real projective $n$-space encloses a volume $V$, it carries a periodic characteristic whose action is at most $\sqrt[n]{n! V}/2$. Moreover, this short characteristic can be chosen so that its projection onto the base manifold is a non-contractible loop.
\end{theorem}

We end the paper by showing that Conjecture~\ref{conj:II} easily implies Mahler's conjectured lower bound for the volume product of centrally symmetric convex bodies: \assertion{if $K \subset \R^n$ is a $0$-symmetric convex body, then $4^n/n! \leq |K| |K^*|$}. This result was partly inspired by the recent preprint of S.~Artstein-Avidan, R.~Karasev, and Y.~Ostrover (\cite{Artstein-Karasev-Ostrover:2013}) who show that the Mahler conjecture is a particular case of C.~Viterbo's volume-capacity conjectured inequality.

\begin{maintheorem}{VII}\label{thm:VII}
    The conjectured sharp isosystolic inequality for reversible optical hypersurfaces in $T^*\RP^n$ implies the Mahler conjecture for centrally symmetric convex bodies in $\R^n$.
\end{maintheorem}

\begin{proof}
    Let $K \subset \R^n$ be the unit ball of a reversible norm $\| \cdot \|$. Without loss of generality we can assume this norm to be a Minkowski norm, the volume being continuous. Consider the flat Finsler metric on $K$ which coincides at each point with the norm $\| \cdot \|$. By identifying the antipodal points in the boundary of $(K,\| \cdot \|)$ we obtain a length space structure on the projective $n$-space and we denote it by $(\RP^n,d_K)$. This length space is not Finsler in general---unless $K$ is an ellipsoid centered at the origin---, and therefore we cannot directly apply the conjectured isosystolic inequality. But we shall see that using this length structure on $\RP^n$ we can define a sequence of smooth Finsler metrics $\{F_k\}$ such that the volume enclosed by their unit co-sphere bundle converges to ${ n! |K| |K^*|}$ and whose homotopical systole---the shortest length of a non-contractible closed curve--- is bigger than $2$. In particular the conjectured isosystolic inequality, applied to these smooth approximations $F_k$, implies when passing to the limit that
    $$
        \frac{4^n}{n!} \leq |K| |K^*| .
    $$
    To construct this sequence of Finsler metrics we strongly perturb the length structure in a small neighborhood of $\partial K$ as follows. First we put the round Riemannian metric on the unit ball $B$ of $\R^n$ that makes it isomorphic to the standard hemisphere of dimension $n$. Then pull back this Riemannian structure by $x \mapsto x/\| x \|$ onto $K \setminus U$, where $U$ is the interior of $\frac{1}{2}K$. We call $F_{aux}$ this metric rescaled in a way that $F_{aux}(x,\cdot) \geq \| \cdot \|$ for all $x \in K \setminus U$. For any integer $k > 2$ let $\phi_k : [0,1] \to [0,1]$ be a smooth, non-decreasing function that is identically zero in $[0,1-{\frac{1}{k}}]$ and identically $1$ in $[1-{\frac{1}{2k}},1]$. Then define the Finsler metric $\tilde{F}_k$ on $K$ by
    $$
        \tilde{F}_k(x,v):=\left(1-\phi_k(\| x \|)\right) \, \| v \|+\phi_k(\| x \|) \, F_{aux}(x, v)
    $$
    for any point $(x,v) \in T K$. Observe that
    \begin{enumerate}
        \item $\| \cdot \| \leq \tilde{F}_k \leq F_{aux}$,
        \item $\tilde{F}_k$ coincides with $\| \cdot \|$ outside the ${\frac{1}{k}}$-neighborhood of $\partial K$ and with $F_{aux}$ inside the ${\frac{1}{2k}}$-neighborhood of $\partial K$.
    \end{enumerate}

    Point $(2)$ implies that $\tilde{F}_k$ is compatible with the identification of antipodal points in $\partial K$ and thus defines a smooth Finsler metric $F_k$ on $\RP^n$. It is straightforward to check that the Holmes-Thompson volume of $F_k$ converges to ${|K| |K^*| / \varepsilon_n}$, and so the volume enclosed by the unit co-sphere bundle $\Sph^*(F_k) \subset T^* \RP^n$ converges to $n! |K||K^*|$.

    It remains to prove that the systole of any $F_k$ is bigger than $2$. Because $\tilde{F}_k \geq \| \cdot \|$, that is enough to prove that the homotopical systole of the length space $(\RP^n,d_K)$ is $2$. It will be more convenient to work on the sphere $S^n$ (the double cover of projective $n$-space) provided with the length metric obtained through its identification with two copies of $K$---each provided with the flat metric induced from the norm $\| \cdot \|$---glued along the boundary.

    Let us denote by $f : S^n \to K$ the natural projection from the two glued copies of $K$ onto $K$. If $\gamma \subset S^n$ is the lift of a non-contractible loop in the projective $n$-space, its image $f(\gamma) \subset K$ is a closed, connected curve that is symmetric about the origin and intersects the boundary of $K$ in at least two antipodal points $x$ and $-x$. Because of the convexity of $K$, the linear projection whose image is the one-dimensional subspace containing these points and whose kernel is parallel to support hyperplanes of $\partial K$ at $x$ and $-x$ is distance (and length) decreasing. It follows that the length of $\gamma$, which equals the length of $f(\gamma)$, is at least $4$. In other words, the length of every non-contractible loop in $(\RP^n, d_K)$ is at least $2$. In fact, the systole of $(\RP^{n},d_{K})$ is equal to $2$ because the length of every diameter of $K$ is equal to $2$ and, once we identify its endpoints, a diameter becomes a non-contractible curve in the projective $n$-space.
\end{proof}

Theorem~\ref{thm:VII} and Ivanov's extension of Pu's inequality to reversible Finsler metrics in~\cite{Ivanov:2002} yield a new proof of Mahler's inequality for centrally symmetric convex bodies on the plane: \assertion{the product of the areas of a $0$-symmetric convex body on the plane and its dual body is greater than or equal to $8$}.

\def\cprime{$'$}
\providecommand{\bysame}{\leavevmode\hbox to3em{\hrulefill}\thinspace}
\providecommand{\MR}{\relax\ifhmode\unskip\space\fi MR }
\providecommand{\MRhref}[2]{%
  \href{http://www.ams.org/mathscinet-getitem?mr=#1}{#2}
}
\providecommand{\href}[2]{#2}


\end{document}